\documentclass[11pt]{article}

\usepackage{ucs}
\usepackage[utf8x]{inputenc}
\usepackage[T1]{fontenc}
\usepackage{color}
\usepackage{caption,subcaption}
\usepackage{graphicx}
\usepackage{amsmath,amsfonts,amssymb}
\usepackage{mathrsfs} 

\definecolor{vertfonce}{rgb}{0.20, 0.46, 0.25}
\definecolor{rougefonce}{rgb}{0.64, 0.09, 0.20}
\newtheorem{proposition}{Proposition}[section]
\usepackage[breaklinks=true,
            colorlinks=true,
            linkcolor=rougefonce,
            citecolor=vertfonce]{hyperref}

\usepackage{autonum} 

\newcommand{\dd}[1]{\ensuremath{\operatorname{d}\!{#1}}}
\newcommand{\RM}{\mathbb{R}}
\newcommand{\ZM}{\mathbb{Z}}
\newcommand{\QM}{\mathbb{Q}}
\newcommand{\NM}{\mathbb{N}}
\newcommand{\CM}{\mathbb{C}}
\newcommand{\h}{\hbar}
\newcommand{\abs}[1]{\left|#1\right|}
\newcommand{\norm}[1]{\left\|#1\right\|}
\newcommand{\demons}[1][$\!\!$]{\noindent\textbf{\demon\ }\textsl{#1}\textbf{.}~}
\newcommand{\cqfd}{\hfill $\square$\par\vspace{1ex}}

\newcommand{\pscal}[2]{\langle{}#1,\,#2\rangle}
\newcommand{\formel}[1]{[\![#1]\!]}
\newcommand{\ad}[1]{{\operatorname{ad}}_{#1}}

\newcommand{\deriv}[2]{\frac{\partial #1}{\partial #2}}
\newcommand{\phy}{\varphi}
\newcommand{\restr}{\upharpoonright}

\newcommand{\theor}{Theorem}
\newcommand{\defin}{Definition}
\newcommand{\lemma}{Lemma}
\newcommand{\remar}{Remark}

\newcommand{\corol}{Corollary}
\newcommand{\propo}{Proposition}
\newcommand{\demon}{Proof}
\newcommand{\probl}{Problem}
\newtheorem{theo}{\theor}[section]
\newtheorem{theo*}{\theor}
\newtheorem{defi}[theo]{\defin}
\newtheorem{defi*}[theo*]{\defin}
\newtheorem{prop}[theo]{\propo}

\newenvironment{rema}
{\par\vspace{1ex}\refstepcounter{theo}%
\noindent\textbf{\remar~\thetheo} }
{~\hfill\mbox{$\triangle$}\par\vspace{1ex}}
\newenvironment{demo}[1][$\!\!$]
{\demons[#1]\ }
{\cqfd}

\title{Semiclassical Birkhoff-Gustavson normal forms and spectral
  asymptotics for nearly resonant Schrödinger operators}

\author{Abdelkader \textsc{Bourebai}\footnote{Department of
    Mathematics, University of Oran1, Ahmed Ben Bella and Laboratory
    of Fundamental and Applied Mathematics of Oran (LMFAO), BP 1524 El Mnaouar, 31000, Oran, Algeria. E-mails: bourebai.abdelkader@edu.univ-oran1.dz and bourebai.abdelkader@yahoo.com}
  \and Kaoutar \textsc{Ghomari}\footnote{ National Polytechnic School
    of Oran, Maurice Audin and Laboratory of Fundamental and Applied (LMFAO)
    Mathematics of Oran, BP 1523 El Mnaouar, 31000, Oran, Algeria. E-mail: kaoutar.ghomari@enp-oran.dz}  \and San \textsc{Vũ
    Ngọc} 
  \footnote{Univ Rennes, CNRS, IRMAR - UMR 6625, F-35000 Rennes, France. E-mail: san.vu-ngoc@univ-rennes1.fr}}

\begin{document}
\maketitle

\begin{abstract}
  The concept of near resonances for harmonic approximations of
  semiclassical Schrödinger operators is introduced and
  explored. Combined with a natural extension of the
  Birkhoff-Gustavson normal form, we obtain formulas for approaching
  the discrete spectrum of such operators which are both accurate and
  easy to implement. We apply the theory to the physically important
  case of the near Fermi (\emph{i.e.} $1:2$) resonance, for which we
  propose explicit expressions and numerical computations.
\end{abstract}
\begin{footnotesize}
\end{footnotesize}
\section{Introduction and motivation}

We are interested in the description of the discrete spectrum of a
semiclassical Schrödinger operator $\hat{P}_\varepsilon$
\begin{equation}
  \hat{P}_\varepsilon (\h) := -\frac{\h^2}{2}\Delta + V_\varepsilon(x),
  \label{equ:schrodinger}
\end{equation}
where $V_\varepsilon$ is a smooth confining potential on $\RM^n$
depending smoothly on a small parameter $\varepsilon\in\RM$, and
$\h>0$ is a small parameter.

More precisely, by ``confining potential'' we mean that there exist
some real value $E_\infty>0$ and a small $\varepsilon_0>0$ such that
for all $E<E_\infty$, the region
$\bigcup_{\abs{\varepsilon}\leq \varepsilon_0}
V_\varepsilon^{-1}((-\infty,E))$ is bounded in $\RM^n$. Moreover, we
also assume that $V_\varepsilon$ grows at most polynomially at
infinity, in the sense of the usual pseudo-differential symbol
classes, uniformly with respect to $\varepsilon$:
\begin{gather}
  \exists m,\forall N\in\NM,\exists C_N,\forall\abs{\varepsilon}\leq
  \varepsilon_0, \forall\alpha\in\NM^n \text{ s.t. } \abs{\alpha}\leq
  N,\\\forall x\in\RM^n \quad \abs{\partial^\alpha V_\varepsilon(x)}
  \leq C_N (1+\norm{x}^2)^{m/2}\,.
\end{gather}
This ensures that $\hat{P}_\varepsilon$ has discrete spectrum in
$(-\infty , E)$.

Let $E_{0,\varepsilon}=\min V_\varepsilon$. We will assume that this
minimum is reached at a unique point and is non-degenerate; hence when
studying the spectrum of $P_\varepsilon$ near $E_{0,\varepsilon}$, it
is natural to consider the harmonic approximation of
$P_\varepsilon$. As we will see in Section~\ref{sec:preparation}, up
to an error of order $\mathcal{O}(\varepsilon^\infty)$, we may assume
that this approximation is smooth in $\varepsilon$. Let us first
consider $\varepsilon=0$. In the harmonic approximation of $P_0$,
which is the quantization of a quadratic Hamiltonian of the form
\begin{equation}
  H_{2}\left( x,\xi \right) =\sum_{j=1}^{n}\frac{\omega_{j}}{2}\left(
    x_{j}^{2}+\xi_{j}^{2}\right),
\end{equation}
two extreme cases may occur. Either the frequencies
$(\omega_1,\dots,\omega_n)$ are independent over $\QM$ (this is the
\emph{non-resonant case}), or they are, up to some common multiple,
all integers. Of course intermediate cases may happen, see
Definition~\ref{definition:re}. In the first case, a well-known result
of Birkhoff~\cite{birkhoff} (following Poincaré) states that the full
symbol $H$ of $P_0$, which is a perturbation of $H_2$, is formally
completely integrable: there are canonical coordinates $(x,\xi)$ and a
smooth map $f$ such that
\begin{equation}
  H(x,\xi) = f(I_1,\dots, I_n) + \mathcal{O}(x,\xi)^\infty,
  \label{equ:BNF-nonresonant}
\end{equation}
where $I_j$ is the action given by
\[
  I_j := \tfrac{1}{2}\left( x_{j}^{2}+\xi_{j}^{2}\right).
\]
Although the Birkhoff idea was soon used by physicists to deal with
quantum Hamiltonians, a rigorous proof of the quantum validity of the
Birkhoff normal form, in the semiclassical limit, is much more recent,
see~\cite{sjostrand-semi}.  This non-resonant case is stable under
perturbations of order $\mathcal{O}(\varepsilon)$ if the quadratic
term $H_2$ is invariant, in the sense that one can write a combined
Birkhoff normal form in all variables $(x,\xi,\varepsilon)$ (and then
a semiclassical Birkhoff normal form in $(x,\xi,\varepsilon,\h)$ will
hold). If one adds diophantine conditions on the frequencies
$\omega_j$, so that they become badly approximated by rationals, then
it is expected that one can accommodate perturbations of the quadratic
term, and even strengthen the result by using KAM stability, similarly
to the case of diophantine tori in~\cite{san-hitrik-sjostrand}. See
also~\cite{lazutkin-book}.

In this paper, we focus on the resonant case, where the situation is
quite different. For simplicity, we will always consider the fully
resonant case, where all frequencies are integers, up to a common
multiple. There is an extension of the Birkhoff normal form for the
resonant case, which was worked out by
Gustavson~\cite{gustavson,eckhardt} (although the general scheme was
already known to Poincaré --- see also Moser's
paper~\cite{moser-new}), where it was shown that, in addition to the
completely integrable normal form~\eqref{equ:BNF-nonresonant}, another
formal series of \emph{resonant terms} has to be considered, which
makes the resulting series much more difficult to analyze (it will be,
generically, non-integrable, see~\cite{duist-112}). The full
semiclassical analysis of resonant harmonic approximations of general
pseudodifferential operators was carried out in~\cite{san-charles}.

When it comes to $\varepsilon$ perturbations, possibly affecting the
harmonic term $H_2$, one may argue that resonant case is not generic:
for most perturbations of a resonant Hamiltonian $H_0$, the perturbed
Hamiltonian $H_\varepsilon$ will be non-resonant. Thus, it is tempting
to claim that, in most physical situations, one can restrict oneself
to the set of non-resonant $\varepsilon$ and thus stick to the
completely integrable normal
form~\eqref{equ:BNF-nonresonant}. However, this normal form is in
general not convergent~\cite{krikorian22}, and, as $\varepsilon\to 0$,
it is expected that the famous appearance of \emph{small denominators}
will make it more and more divergent, hindering the effectiveness of
the approximation (unless the full Hamiltonian is known to be
integrable and analytic, see~\cite{zung-birkhoff}).

This phenomenon of near resonances has been recognized as crucial in
molecular spectroscopy. Thus, for the study of the dynamics of highly
excited vibrational states, Joyeux shows in his article
\cite{joyeux-birkhoff} that the resonant Birkhoff-Gustavson procedure
can yield more accurate results than the standard Birkhoff procedure,
even when the Hamiltonian is \emph{not} resonant. He considers the HCP
molecule, called phosphatine, where the calculation of the fundamental
frequencies leads to $ \omega_{1}\approx 1256 $ (C-P stretch),
$\omega_{2}\approx 650 $ (bend), and
$\omega_{3}\approx 3479 \textup{cm}^{-1} $ (C-H stretch).  The
Hamiltonian obtained by the non-resonant Birkhoff normal form (also
called the Dunham expansion) is the formal series
\[
  H_{D}=\sum_{i}\omega_{i}I_{i}+\sum_{i\leqslant j}
  x_{ij}I_{i}J_{j}+\sum_{i\leqslant j\leqslant k}
  y_{ijk}I_{i}J_{j}I_{k}+ \dots
\]
Joyeux computed levels of HCP up to $ 22000\textup{cm}^{-1}$ above the
bottom of the well by truncating the series at various orders, and
compared the results to the exact quantum levels of HCP relative to
the ground state.  Using $\ell^1$, $\ell^2$ and $\ell^\infty$ norms to
estimate the discrepancy with respect to the exact quantum
computation, he observed a rapid divergence, see \cite[Table
1]{joyeux-birkhoff}, which limits the interest of this Hamiltonian to
at most $ 4th $ or $ 5th $ order, which is not satisfactory.  The
article therefore concludes that the Dunham expansion is very far for
being sufficient for whatever quantitative purpose.

Given that the relation between both fundamental frequencies
$ \omega_{1} $ and $ \omega_{2} $, namely
$ \omega_{1}-2\omega_{2} \approx -44\textup{cm}^{-1}$ is a \emph{near
  resonance} relation, Joyeux proposed, as a next step, to take this
resonance into account in the formal expansion, and to compare once
again the results of the calculations of the energy levels of the HCP
molecule. For this, he considered the following Hamiltonian:
\[
  H=H_{D}+\sum_{m\geqslant 1}H_{F}^{(m)},
\]
\[
  H_{F}^{(m)}= 2I_{1}^{ \frac{m}{2}} I_{2}^{m} \cos(m \varphi_{1} - 2m
  \varphi_{2}) \times \left( k^{(m)}+ \sum_{i}k_{i}^{(m)}I_{i}+
    \sum_{i\leqslant j} k_{ij}^{(m)}I_{i}J_{j}+ \dots\right)\,,
\]
which he calls the Fermi resonance Hamiltonian (in Physics or
Chemistry literature, the $2:1$ resonance is traditionally called the
Fermi resonance, see Section~\ref{sec:b.g.n.f-near-fermi}). The terms
in $H_F$ are precisely those given by the Birkhoff-Gustavson procedure
in case of an exact $2:1$ resonance. It is then noticed that the
results of the computations of energy levels are much more accurate
with this modified Hamiltonian.

From this study comes the motivation and main goal of our work, which
is a description of the spectrum of semiclassical Schrödinger
operators $ \hat{P}_\varepsilon (\h)$ given
in~\eqref{equ:schrodinger}, for which the harmonic frequencies are
close to resonance.  For this purpose, we build on the
paper~\cite{san-charles} which gives precise semiclassical asymptotics
for an exact, full resonance, by restricting the Hamiltonian to the
eigenspaces of the resonant harmonic oscillator (whose dimensions tend
to infinity in the semiclassical limit).
Our main results are organized as follows:
 
In Section~\ref{sec:preparation}, we consider the first step of this
work, which consists in transforming the initial Schrödinger
Hamiltonian $ \hat{P}_{\varepsilon}$ into a perturbation of the
harmonic oscillator (Proposition~\ref{prop:preparation}). For this, we
need a diagonalization result of symmetric real matrices depending
smoothly on small parameter $ \varepsilon $ (Theorem~\ref{theo:symm}).
 
Section~\ref{sec:BGNF} has a double goal. On the one hand, we prove in
Theorem~\ref{theo:BGNF} that the Birkhoff-Gustavson normal form (BGNF)
theorem can be extended to handle Schrödinger operators
$ \hat{P}_{\varepsilon}$ depending on the parameter $ \varepsilon $
(where the harmonic frequencies also depend on the parameter
$ \varepsilon $); on the other hand, we give an explicit construction
of the BGNF in the near Fermi resonance, in
Theorems~\ref{theo:mu-reso12} and~\ref{theo:mu_explicit}.

In Section~\ref{sec:spectr-analys-fermi}, we give the exact matrix
representation of the ``polyads'' generated by the first non-trivial
Birkhoff correction of the Fermi resonance, \emph{i.e.} the
restriction of the quantum BGNF to the various eigenspaces of the
resonant harmonic oscillator (Theorem~\ref{theo:matrix}).

Finally, in Section~\ref{sec:numerics}, we show how the theoretical
study can lead to numerical schemes, and we propose, in the case of
the Fermi resonance, a detailed numerical illustration of our results,
by comparing the ``exact'' quantum spectrum of $\hat{P}_\varepsilon$
with the eigenvalues obtained via the $\varepsilon$-Birkhoff-Gustavson
procedure, at order 3.

\section{Preparation: smooth diagonalization}
\label{sec:preparation}

In order to obtain the harmonic approximation of the Schrödinger
operator~\eqref{equ:schrodinger}, we need to diagonalize the Hessian
of $V_\varepsilon$ at the critical point, in a smooth way. Because our
aim is to deal with resonant eigenvalues, we cannot assume that
eigenvalues are simple, and we will use the following general result,
which is elementary but apparently not often found in the literature
(we could not locate a reference).

\begin{theo}
  \label{theo:symm}
  Let $A(\varepsilon)$ be a family of $n\times n$ real symmetric,
  respectively hermitian, matrices depending in a smooth (\emph{ie.}
  $C^{\infty}$) way on a small parameter $\varepsilon$.  Then there
  exists a smooth family of orthogonal, respectively unitary, matrices
  $U(\varepsilon)$ and smooth functions
  $\varepsilon\mapsto \lambda_{j,\varepsilon}\in\RM$, $j=1,\dots,n$
  such that
  \begin{equation}
    U^{-1}(\varepsilon)A(\varepsilon)U(\varepsilon) = 
    \textup{diag}(\lambda_{1,\varepsilon},\dots,\lambda_{n,\varepsilon}) + 
    \mathcal{O}(\varepsilon^{\infty}).
  \end{equation}
\end{theo}
\begin{rema}
  This result means that, if we accept to replace the true eigenvalues
  by approximate eigenvalues which are close to the exact one to any
  order in $\varepsilon$, then we may smoothly diagonalize the family
  $A(\varepsilon)$. What may comes as a surprise when one first
  encounters this kind of result is that, in general, it is \emph{not}
  possible to smoothly diagonalize a smooth family of symmetric
  matrices in an exact way (unless all eigenvalues are simple). See
  for instance Example 5.3, section II-5 of the Kato
  book~\cite{Kato}. On the other hand, positive results are available
  when the family is analytic, or if one only requires $C^1$
  regularity for the eigenvalues, a result due to
  Rellich~\cite{rellich}, see also~\cite[Theorem 6.8 Section
  II-6]{Kato}.
\end{rema}
\begin{demo}
  Let us treat the real symmetric case; the Hermitian case will be
  completely analogous.  We reason by induction on the size of the
  matrix. The result obviously holds when $n=1$ (without the
  $\mathcal{O}(\varepsilon^\infty)$ error term).  We may assume that
  $A(0)$ is diagonal; let $\mu_1<\dots<\mu_\ell$ be its eigenvalues,
  with multiplicities $d_1,\dots, d_\ell$. By the min-max formula, the
  eigenvalues of $A(\varepsilon)$ are continuous in
  $\varepsilon$. Thus, there exists $\rho>0$ such that the spectrum of
  $A(\varepsilon)$ is contained in
  $\bigcup_{j=1}^\ell B(\mu_{j},\rho) $, and each ball
  $ B_{j} := B(\mu_{j},\rho) $ contains $ d_{j} $ eigenvalues (counted
  with multiplicity). This shows that the spectral projector on the
  generalized eigenspaces,
  \[
    P_{B_{j}}(\varepsilon) := -\frac{1}{2\pi i}\int_{\partial
      B_{j}}(A(\varepsilon)-z)^{-1}dz,
  \]
  is $ C^{\infty} $ in $ \varepsilon $. It is now easy to find an
  orthonormal basis of the generalized eigenspace
  $E_{j}(\varepsilon):=P_{B_{j}}(\RM^n)$ that depends smoothly on
  $\varepsilon$. For instance, one can take a basis $\mathcal{B}(0)$
  of $E_{j}(0)$; for $\varepsilon$ small enough, the projection
  $P_{B_j}(\varepsilon)(\mathcal{B}(0))$ is a basis of
  $E_j(\varepsilon)$, which we may smoothly orthonormalize by the
  Gram-Schmidt algorithm. In this way we obtain a smooth bloc-diagonal
  decomposition: there exists a smooth unitary matrix
  $ V(\varepsilon) $ such that
  \[
    V^{-1}(\varepsilon) A(\varepsilon) V(\varepsilon) = \left(
      \begin{array}{ccc}
        A_{1}(\varepsilon) & 0 & 0\\
        0 & \ddots & 0\\
        0 & 0 & A_{\ell}(\varepsilon)
      \end{array}
    \right) \;.
  \]
  where $ A_j(\varepsilon) = \mu_{j}I+\mathcal{O}(\varepsilon)$ is a
  real symmetric matrix. If $\ell>1$ we obtain the result by
  induction.  Hence it only remains to consider the case of a unique
  generalized eigenspace of dimension $n$: we have
  \[
    A(\varepsilon) = \mu_1 I + \mathcal{O}(\varepsilon).
  \]
  Because the remainder $\mathcal{O}(\varepsilon)$ is smooth, we can
  write it as $\varepsilon B_1(\varepsilon)$, where $B_1(\varepsilon)$
  is smooth (and real symmetric). Therefore the question is reduced to
  diagonalizing $B_1(\varepsilon)$, and we may repeat the
  procedure. There are finally only two possibilities:
  \begin{enumerate}
  \item Either there exists $N>0$ such that, after the $N$-th
    iteration, the remainder $B_N(0)$ possesses more than one
    generalized eigenspace. Then we may split them and obtain the
    result as above;
  \item or for all $N>0$, $B_N(0)$ has only one generalized
    eigenspace.
  \end{enumerate}
  In the second case, we have real constants $c_0=\mu_1,c_1,c_2,\dots$
  such that, for any $N\geq 0$,
  \[
    A(\varepsilon)=c_{0} I+c_{1}\varepsilon
    I+\dots+c_{N}\varepsilon^{N}I + \mathcal{O}(\varepsilon^{N+1}).
  \]
  By the Borel Lemma, there exists a $ C^{\infty} $ function
  $ \lambda(\varepsilon) $ whose Taylor series at $\varepsilon=0$ is
  \[
    \lambda(\varepsilon) \sim
    c_{0}+c_{1}\varepsilon+c_{2}\varepsilon^{2}+\dots+c_{N}\varepsilon^{N}+\dots
  \]
  which means
  $ \forall N, A(\varepsilon) - \lambda(\varepsilon)I =
  \mathcal{O}(\varepsilon^{N})$.  Thus,
  $A(\varepsilon) = \lambda(\varepsilon)I +
  \mathcal{O}(\varepsilon^{\infty})$, which gives the result.
\end{demo}

\section{Birkhoff-Gustavson normal form in near resonance}
\label{sec:BGNF}

In this section, we shall discuss the Birkhoff normal form procedure
for Schrödinger operators $ \hat{P}_\varepsilon $ which depend on
small parameters $ {\h} >0 $ and $ \varepsilon $, and we will then
apply the general ideas to the near Fermi resonance.

\subsection{$\varepsilon$-Birkhoff-Gustavson normal form theorem}

On $ L^2(\RM^n) $ consider the Schrödinger operator
\[
  \hat{P_\varepsilon} = -\tfrac{\h^2}{2}\Delta_x + V_{\varepsilon}(x),
\]
where $ \Delta =\Delta_x$ is the $ n $ dimensional Laplacian and
$ V_{\varepsilon} $ is a smooth real potential depending smoothly on
$ \varepsilon $.  We wish to perform a local (and microlocal) analysis
near the origin $x= 0 $. To this effect, we assume that for
$ \varepsilon=0 $, the potential $ V_{0} $ has a non-degenerate
minimum at the origin:
\[
  V_{0}(0)=0,V_{0}'(0)=0,V_{0}''(0)>0.
\]
Using the implicit function theorem to
$ F(\varepsilon, x)=V_{\varepsilon}'(x) $, we obtain a smooth map
$\varepsilon\mapsto x_{\varepsilon}$ near the origin such that $x_0=0$
and $V_{\varepsilon}'(x_{\varepsilon})=0$, for $ \varepsilon>0 $ small
enough. So, there exists $\varepsilon_0>0$ such that for all
$\abs{\varepsilon}\leq \varepsilon_0$, the point $ x_{\varepsilon} $
is a non-degenerate minimum for $ V_{\varepsilon} $.  Using the
translation $x\mapsto \tilde x =x-x_{\varepsilon} $, which yields a
unitary map $\tau_{\varepsilon}$ on $L^2(\RM^n)$, given by
$\tau_\varepsilon f(\tilde x) = f(\tilde x+x_\varepsilon)$,
$\hat{P}_\varepsilon $ is transformed into
\[
  \tau_{\varepsilon}\hat{P}_\varepsilon \tau_{\varepsilon}^{-1} =
  -\tfrac{\h^2}{2}\Delta_{\tilde x} + W_{\varepsilon}(\tilde x),
\]
where
$ W_{\varepsilon}(\tilde x)=V_{\varepsilon}(\tilde x+x_\varepsilon)
$. We have
\[
  W_{\varepsilon}(0)=V_{\varepsilon}(x_\varepsilon), \quad
  W_{\varepsilon}'(0)=V_{\varepsilon}'(x_{\varepsilon})=0\,,
\]
and the symmetric matrix $ W_{\varepsilon}''(0)$ is positive definite,
for $ \varepsilon>0 $ sufficiently small.

Using theorem \ref{theo:symm}, we can smoothly diagonalize
$ W_{\varepsilon}''(0) $ modulo $ \mathcal{O}(\varepsilon^{\infty}) $,
via a change of variables $y= U_\varepsilon^*\tilde x$, with
$U_\varepsilon\in O(n)$.  Let
$ \left(\omega_{1,\varepsilon}^2,\dots,
  \omega_{n,\varepsilon}^2\right) $ be its approximate eigenvalues:
they are positive and depend in a $ C^\infty $ way on $ \varepsilon $.

The change of variables $U_\varepsilon$ induces a unitary map
$g_\varepsilon$ on $L^2(\RM^n)$ given by
$g_\varepsilon f(y) = f(U_\varepsilon y)$. Since $U_\varepsilon$ is
orthogonal,
$g_\varepsilon \Delta_{\tilde x} g_\varepsilon^{-1} = \Delta_y$ and
hence the transformed operator is
\[
  g_\varepsilon \tau_{\varepsilon} \hat{P}_\varepsilon
  \tau_{\varepsilon}^{-1}g_\varepsilon^{-1} =
  -\tfrac{\h^2}{2}\Delta_{y} + W_{\varepsilon}(U_\varepsilon y),
\]
Since
$W_\varepsilon(\tilde x) = \frac{1}{2}\pscal{W_{\varepsilon}''(0)
  \tilde x}{\tilde x} + \mathcal{O}(\tilde x^3)$, we obtain
\begin{equation}
  W_{\varepsilon}(U_\varepsilon y) = W_\varepsilon(0) +
  \frac{1}{2} \sum_{j= 1}^n \omega_{j,\varepsilon}^2 y_j^2 +
  \mathcal{O}(\varepsilon^\infty\norm{y}^2) + \mathcal{O}(y^3)\,.
  \label{equ:Wepsilon}
\end{equation}

The rescaling
$ y_{j}\mapsto \tilde y_j :=\sqrt{\omega_{j,\varepsilon}}y_{j} $,
giving rise to the unitary map
$\Lambda_\varepsilon:f\mapsto
(\tilde{y}\mapsto\frac{1}{(\omega_{1,\varepsilon}\cdots
  \omega_{n,\varepsilon})^{n/2}}f(\frac{\tilde{y}_1}{\sqrt{\omega_{1,\varepsilon}}},
\dots, \frac{\tilde{y}_n}{\sqrt{\omega_{n,\varepsilon}}}))$ transforms
$ \hat{P_\varepsilon} $ into a perturbation of the harmonic oscillator
$ \hat{H}_{2,\varepsilon} $:
\begin{equation}
  \label{eq:reste}
  \Lambda_\varepsilon g_\varepsilon \tau_{\varepsilon}\hat{P}_{\varepsilon} \tau_{\varepsilon}^{-1}g_\varepsilon^{-1}\Lambda_\varepsilon^{-1} =
  W_{\varepsilon}(0)+\hat{H}_{2,\varepsilon}+R_{\varepsilon}(\tilde{y}),
\end{equation}
where,
\[
  \hat{H}_{2,\varepsilon}=
  \sum_{j=1}^{n}\frac{\omega_{j,\varepsilon}}{2}\left(
    -\h^{2}\frac{{\partial}^2}{{\partial \tilde{y}_{j}^2}} +
    \tilde{y}_{j}^{2}\right)
  +\mathcal{O}(\varepsilon^\infty\norm{\tilde{y}}^2),
\]
and $ R_{\varepsilon}(\tilde y)$ is a smooth function of order
$ \mathcal{O}( \tilde y ^3) $ at the origin, uniformly with respect to
$\varepsilon$.

From now on, to simplify notation, we switch back to $x$, and assume
that
\[
  \hat{P}_{\varepsilon}=
  W_{\varepsilon}(0)+\hat{H}_{2,\varepsilon}+R_{\varepsilon}(x),
\]
with
\[
  \hat{H}_{2,\varepsilon}=
  \sum_{j=1}^{n}\frac{\omega_{j,\varepsilon}}{2}\left(
    -\h^{2}\frac{{\partial}^2}{{\partial x_{j}^2}} + x_{j}^{2}\right)
  +\mathcal{O}(\varepsilon^\infty\norm{x^2})\,.
\]

Let
\[
  \hat{H}_{2,0}=
  \sum_{j=1}^{n}\frac{\omega_{j,0}}{2}\left(-\h^{2}\frac{{\partial}^2}{{\partial
        x_{j}^2}}+x_{j}^{2}\right),
\]
The eigenvalues $\omega_{j,\varepsilon}$ being smooth in
$\varepsilon$, we have
\begin{equation}
  \label{eq:frequences}
  \omega_{j,\varepsilon}=\omega_{j,0}+\varepsilon \omega_{j,1}+{\varepsilon}^2 {\tilde {\omega}}_{j,2}\left(\varepsilon \right)\,.
\end{equation}

Therefore,
\begin{equation}
  \hat{P}_{\varepsilon}
  = W_{\varepsilon}(0)+\hat{H}_{2,0} +
  \varepsilon \hat{L}_{2}+ {\varepsilon}^{2}\hat M_{2,\varepsilon} +
  R_{\varepsilon}(x)+\mathcal{O}(\varepsilon^\infty\norm{x^2}),
\end{equation}
where
\begin{equation}
  \hat{L}_{2} = \sum_{j=1}^{n} \frac{\omega_{j,1}}{2} \left(
    -\h^{2}\frac{{\partial}^2}{{\partial x_{j}^2}}+x_{j}^{2}\right) \
  \text{and}\ \ \hat M_{2,\varepsilon}=\sum_{j=1}^{n}
  \frac{{\tilde{\omega}}_{j,2}(\varepsilon)
  }{2}\left(-\h^{2}\frac{{\partial}^2}{{\partial x_{j}^2}}+x_{j}^{2}
  \right),\label{equ:hat-L2-M2}
\end{equation}
Because the symbol of $ \hat{L}_{2} $ is quadratic, it cannot be
reduced by the usual Birkhoff-Gustavson procedure; in order to deal
with this issue, we have to add $\varepsilon$ to the set of formal
variables, so that $\varepsilon \hat{L}_{2}$ becomes of order 3.

Thus, we introduce the space
\[
  \mathcal{E}=\CM \formel{x,\xi,\varepsilon,\h} =\CM \formel{
    x_{1},\dots,x_{n},\xi_{1},\dots,\xi_{n},\varepsilon,\h} ,
\]
of formal power series of $ (2n+2) $ variables with complex
coefficients, where the degree of the monomial
$ x^{\alpha}\xi^{\beta}\varepsilon^{m}\h^{\ell} $ is defined to be
$ \vert\alpha\vert+\vert\beta\vert+m+2\ell $,
$ \alpha,\beta\in \NM^n,\ell,m \in \NM $.

Let $ \mathcal{D}_{N} $ be the finite dimensional vector space spanned
by monomials of degree $ N $ and $\mathcal{O }_N$ the subspace of
$ \mathcal{E} $ consisting of formal series whose coefficients of
degree $ < N $ vanish, then
$\left(\mathcal{O}_{N} \right)_{N\in \NM} $ is a filtration:
\[
  \mathcal{E}\supset\mathcal{O}_{0}\supset\mathcal{O}_{1}\supset\cdots,
  \qquad \bigcap_{N}\mathcal{O}_{N}=\left\lbrace 0 \right\rbrace\,.
\]
This filtration will be used for all formal convergences in this
section. We shall also need to discuss the degree in $(x,\xi,\h)$
only; to this aim, we denote by $\mathcal{O}_N(x,\xi,\h)$ the subspace
of $\mathcal{E}$ spanned by monomials
$ x^{\alpha}\xi^{\beta}\varepsilon^{m}\h^{\ell} $ such that
$\abs{\alpha}+ \abs{\beta} + 2\ell\geq N$.

Let $ d,m \in{\RM} $, the symbol class $S^{d}(m)$ is the set of the
smooth functions
$
a_{\varepsilon}\left(.,\h\right):{\RM}^{2n}\times]0,1]\rightarrow\CM$
such that, for all $ \alpha \in \NM^{n} $,
\[
  \vert\partial_{(x,y)}^{\alpha}a_{\varepsilon} (x,\xi,\h) \vert\leq
  C_{\alpha}\h^d \left(1+\vert x\vert ^{2}+\vert \xi \vert
    ^{2}\right)^{\frac{m}{2}}
\]
for some constant $ C_{\alpha}>0$, uniformly in $\h\in (0,\h_0]$ and
$ \varepsilon\in[-\varepsilon_0,\varepsilon_0]$, where $\h_0$ and
$\varepsilon_0$ are small enough. $S^{d}(m)$ is called the space of
symbols of order $ d $ and degree $ m $.

For $ a_\varepsilon\in S^{d}(m)$ and $ u\in C^{\infty}_0(\RM^{2n}) $,
the Weyl quantization of $a_\varepsilon$ acting on $u$ is given by the
oscillatory integral
\begin{equation}
  (Op_{\h}^{W}(a_{\varepsilon})u)(x)=\frac{1}{(2 \pi \h)^{n}}\int_{{\RM}^{2n}}
  e^{\frac{i}{\h}\langle x-y,\xi\rangle}a_{\varepsilon} \left(\tfrac{x+y}{2},\xi ,\h \right)u(y) \vert dyd\xi \vert.
\end{equation}
In general, $ Op_{\h}^{W}(a_\varepsilon) $ is an unbounded linear
operator on $ L^2(\RM^n) $, and $a_\varepsilon$ is called its Weyl
symbol. For example, the Weyl symbol of the harmonic oscillator
$ \hat{H}_{2,0} $ is the polynomial
\begin{equation}
  H_{2,0}=
  \sum_{j=1}^{n}\frac{\omega_{j,0}}{2}(\xi_j^{2}+x_{j}^{2}) \,,
  \label{eq:H_20}
\end{equation}
and the Weyl symbol of the operator of multiplication by a function
$f(x)$ is the function $f$ itself or, more precisely, the function
$(x,\xi)\mapsto f(x)$.

In this article, we use the Weyl bracket $ \left[f,g\right]_{W}$
defined on $\mathcal{E}$ by the formal Taylor series of the Weyl
symbol of the commutator
\[
  Op_{\h}^{W} \tilde f \circ Op_{\h}^{W} \tilde g - Op_{\h}^{W} \tilde
  g \circ Op_{\h}^{W} \tilde f,
\]
where $\tilde f$ and $\tilde g$ are smooth symbols whose formal Taylor
series is equal to $f$ and $g$, respectively. From now on, when
$f=f(x,\xi,\varepsilon,\h)$ is a smooth function, we shall allow us to
write $f\in\mathcal{E}$ to signify that its Taylor series at the
origin belongs to $\mathcal{E}$.

The filtration $\mathcal{E}$ has a nice behaviour with respect to the
Weyl bracket: if $ N_{1}+N_{2}\geqslant 2$, then,
$ \h ^{-1}[\mathcal{O}_{N_{1}},\mathcal{O}_{N_{2}}]_{W} \subset
\mathcal{O}_{N_{1}+N_{2}-2} $.

The Weyl symbol $ H_{2,\varepsilon} $ of $ \hat{H}_{2,\varepsilon} $
is:
\begin{align}\label{equ:L2-M2}
  H_{2,\varepsilon}
  &= H_{2,0}+\varepsilon L_{2}+\varepsilon^2 M_{2,\varepsilon}
    +\mathcal{O}(\varepsilon^\infty\norm{x}^2) \\
  &= \sum_{j=1}^{n}\frac{\omega_{j,0}}{2}\left(\xi_j^{2}+x_{j}^{2}\right)
    + \varepsilon\sum_{j=1}^{n}\frac{\omega_{j,1}}{2}\left(\xi_{j}^{2}+x_{j}^{2} \right)+{\varepsilon}^{2}M_{2,\varepsilon}
    + \mathcal{O}(\varepsilon^\infty\norm{x}^2),
\end{align}
with
\[
  L_{2}=\sum_{j=1}^{n}\frac{\omega_{j,1}}{2}\left(\xi_{j}^{2}+x_{j}^{2}\right)\,,
\]
where $M_{2,\varepsilon}$ is the Weyl symbol of
$\hat M_{2,\varepsilon}$ and $L_2$ is the Weyl symbol of $\hat L_2$.

To summarize, we have proven the following.
\begin{prop}
  \label{prop:preparation}
  Let $ \hat{P}_{\varepsilon} $ be the Schrodinger operator defined in
  \eqref{equ:schrodinger}. Then there exists an explicit unitary
  operator $U_\varepsilon$ on $L^2(\RM^n)$
  ($U_\varepsilon= \Lambda_\varepsilon g_\varepsilon
  \tau_{\varepsilon}$ is composed of a translation, a unitary
  transform, and a scaling, all in the position variable $x$) such
  that the differential operator
  $U_\varepsilon \hat{P}_{\varepsilon} U_\varepsilon^* $ has the
  following Weyl symbol:
  \begin{align}
    \sigma_W(U_\varepsilon \hat{P}_{\varepsilon} U_\varepsilon^* )
    & = \min_{\RM^n} V_\varepsilon + H_{2,0}+\varepsilon L_{2}+L_{\varepsilon}',\\
    & =  \min_{\RM^n} V_\varepsilon + H_{2,0}+L_{\varepsilon},
  \end{align}
  where $H_{2,0}$ and $L_2$ are defined in~\eqref{equ:L2-M2},
  $ L_{2}\in \mathcal{D}_{2}$,
  $L_{\varepsilon}' = \varepsilon^2M_{2,\varepsilon} +
  R_{\varepsilon}(x)+\mathcal{O}(\varepsilon^ \infty\norm{x}^2)\in
  \mathcal{O}_{3} $,
  $L_{\varepsilon}= \varepsilon L_{2}+L_{\varepsilon}'\in
  \mathcal{O}_{3} $,$M_{2,\varepsilon}$ is the Weyl symbol of
  $\hat{M}_{2,\varepsilon}$ defined in~\eqref{equ:hat-L2-M2}, and and
  $ R_{\varepsilon}(x)$ is a smooth function of order
  $ \mathcal{O}(x^3) $ at the origin, uniformly with respect
  to $\varepsilon$.
\end{prop}
    
Let $A\in \mathcal{E}$, we define the adjoint operator
$\ad{A}(P):=[A,P]_{W}$, $P\in \mathcal{E}$. The crucial properties of
$ \ad{H_{2,0}} $ are given in the next proposition, for more details
see~\cite {san-charles}.
\begin{prop}
  \label{prop:freq}
  \begin{enumerate}
  \item Let $\omega_0=(\omega_{1,0},\dots,\omega_{n,0})$, and let
    $z_j=x_j+ i \xi_j$; then $ \ad{H_{2,0}} $ is diagonal on the
    family
    $\lbrace z^{\alpha}\bar{z}^{\beta},z\in{\CM}^{n},\alpha,\beta\in
    {\NM}^{n}\rbrace $ of
    $ \CM\formel{z,\bar{z}} = \CM\formel{x,\xi} $ and
    \[
      \ad{H_{2,0}}(z^{\alpha}\bar{z}^{\beta}) = \h \langle
      \beta-\alpha,\omega_0\rangle z^{\alpha}\bar{z}^{\beta} \, .
    \]
  \item We have
    \begin{equation}
      \mathcal{D}_{N}= \ker \left(\h^ {-1} \ad{H_{2,0}} \right) \oplus \textup{Im}
      \left(\h^ {-1} \ad{H_{2,0}} \right).\label{equ:admissible}
    \end{equation}
  \end{enumerate}
\end{prop}
We are now ready to state the formal $\varepsilon$-Birkhoff-Gustavson
normal form, which is a natural extension of the usual
Birkhoff-Gustavson case. Note that we prove here the result for a
general perturbation $L_\varepsilon\in\mathcal{O}_3$, which does not
need to be equal to the particular form obtained above when one
reduces a Schrödinger operator.
\begin{theo}
  \label{theo:BGNF} Let $H_{2,0}$ be as in~\eqref{eq:H_20}, and let
  $L_{\varepsilon}\in \mathcal{O}_{3}$, then there exist
  $A_{\varepsilon}\in \mathcal{O}_{3}$ and
  $K_{\varepsilon}\in \mathcal{O}_{3}$ such that,
  \begin{equation}
    \label{eq:B.G.N.F}
    \begin{split}
      e^{i\h^{-1}\ad{A_{\varepsilon}}} \left(H_{2,0}+L_{\varepsilon} \right)&=H_{2,0}+K_{\varepsilon},\\
      \left[H_{2,0},K_{\varepsilon} \right]_W=0.
    \end{split}
  \end{equation}
  Moreover, the following properties hold.

  \begin{enumerate}
  \item $K_\varepsilon$ is unique, and $A_{\varepsilon}$ is unique
    modulo $\ker \ad{H_{2,0}}$.
  \item if $H_{2,0}$ and $L_{\varepsilon}$ have real coefficients
    (\emph{i.e.} belong to $\RM \formel{x,\xi,\varepsilon,\h}$) then
    $K_\varepsilon\in \RM \formel{x,\xi,\varepsilon,\h}$, and
    $A_{\varepsilon}$ can be chosen to have real coefficients as well.
  \item If
    $L_\varepsilon\in \varepsilon^m \CM \formel{x,\xi,\varepsilon,\h}
    $, for some $m\in\NM$, then
    $K_\varepsilon\in \varepsilon^m \CM \formel{x,\xi,\varepsilon,\h}
    $ as well.
  \item If
    $L_\varepsilon\in \h^\ell \CM \formel{x,\xi,\varepsilon,\h} $, for
    some $\ell\in\NM$, then
    $K_\varepsilon\in \h^\ell \CM \formel{x,\xi,\varepsilon,\h} $ as
    well.
  \item If $L_\varepsilon\in\mathcal{O}_N(x,\xi,\h)$ for some
    $N\geq 0$, then $K_\varepsilon\in\mathcal{O}_N(x,\xi,\h)$ as well.
  \end{enumerate}
\end{theo}
Notice that the sum :
\begin{equation}
  e^{i\h^{-1}\ad{A_{\varepsilon}}}\left(H_{2,0}+L_{\varepsilon} \right)=\sum_{\ell=0}^{\infty}\frac{1}{\ell!} \left(\tfrac{i}{\h} \ad{A_{\varepsilon}} \right)^{\ell}\left(H_{2,0}+L_{\varepsilon} \right),
\end{equation}
is usually not convergent in the analytic sense, even if
$ P_{\varepsilon} $ is analytic, but it is always convergent in the
formal topology of $ \mathcal{E} $, because the map
$ B\mapsto\frac{i}{\h} \ad{A}(B)=\frac{i}{\h} \left[A,B \right]_{W}$
sends $\mathcal{O}_{N}$ into $\mathcal{O}_{N+1}$.

\medskip

\begin{demo}[of the Theorem]
  We construct $A_{\varepsilon} $ and $ K_{\varepsilon} $ inductively,
  by successive approximations with respect to the filtration
  $ \mathcal{E}.$

  Let $ N \geq 1 $, we suppose that there exists
  $ A_{N, {\varepsilon}} \in \mathcal{O}_{3} $ and
  $ K_{N, {\varepsilon}} \in \mathcal{O}_{3} $ such:
  \begin{equation}\label{equ:inductionK}
    e^{i\h^{-1} \ad{A_{N,{\varepsilon}}}}\left(H_{2,0}+L_{\varepsilon} \right)=H_{2,0}+K_{3, \varepsilon}+ \dots +K_{N+1, \varepsilon}+R_{N+2, \varepsilon}+ \mathcal{O}_{N+3}
  \end{equation}
  where $ K_{j, \varepsilon} \in \mathcal{D}_{j}$ commutes with
  $ H_{2,0} $ and $ R_{N+2, \varepsilon} \in \mathcal{D}_{N+2}. $ We
  look for $ A_{\varepsilon}' \in \mathcal{D}_{N+2} $, and
  $ K_{N+2, \varepsilon} \in \mathcal{D}_{N+2} $ that commutes with
  $ H_{2,0} $, such that :
  \begin{equation}
    e^{i\h^{-1} \ad{A_{N,{\varepsilon}}+A_{\varepsilon}'}}\left(H_{2,0}+L_{\varepsilon} \right)=H_{2,0}+K_{3, \varepsilon}+ \dots +K_{N+1, \varepsilon}+K_{N+2, \varepsilon}+ \mathcal{O}_{N+3}\,.
  \end{equation}
  This is equivalent to
  \begin{equation}
    R_{N+2, \varepsilon}-\frac{i}{\h}\left[ H_{2,0},A_{\varepsilon}' \right]_{W}=K_{N+2, \varepsilon}.
  \end{equation}
  Since $ H_{2,0} $ satisfies~\eqref{equ:admissible}, we have
  $ R_{N+2, \varepsilon}=R_{N+2, \varepsilon}'+\frac{i}{\h}\left[
    H_{2,0},R_{N+2\varepsilon}'' \right]_{W} $, where
  $ R_{N+2, \varepsilon}' \in \mathcal{D}_{N+2}$ commutes with
  $ H_{2,0} $ (and is unique) and
  $ R_{N+2,\varepsilon}'' \in \mathcal{D}_{N+2} $, unique modulo
  $\ker \ad{H_{2,0}}$.  Thus, we must (and can) choose
  $ K_{N+2, \varepsilon}= R_{N+2, \varepsilon}'$ and, we may choose
  $ A_{\varepsilon}'= R_{N+2,\varepsilon}'' $. This shows
  that~\eqref{equ:inductionK} holds at order $N+1$,
  proving~\eqref{eq:B.G.N.F}, and the first property.

  The second property follows from Moyal's formula
  \[
    i [ f,g ]_W (x, \xi , \h ) = 2\sin \Bigl( \frac{\h}{2} \square
    \Bigr)\bigl( f ( x_1, \xi_1, \h ) g ( x_2, \xi_2, \h ) \bigr)
    \Bigr|_{x=x _1 = x_2,\atop \xi = \xi_1 = \xi_2}
  \]
  where
  \[
    \square = \sum_{j=1}^n \partial_{\xi^j_1} \partial_{x^j_2} -
    \partial_{x^j_1} \partial_{\xi^j_2}\,.
  \]
  Properties 3 and 4 are due to the fact that $\varepsilon$ and $\h$
  are central elements in $\mathcal{E}$ (they commute with
  everything). The last properties holds because $\ad{H_{2,0}}$
  preserves the $(x,\xi,\h)$-order.
  
\end{demo}

\subsection{Normalizing the near Fermi resonance}
\label{sec:b.g.n.f-near-fermi}

A vector of $n$ frequencies $(\omega_1,\dots,\omega_n)$ is called
\emph{resonant} if the coefficients $\omega_1,\dots,\omega_n$ are
linearly dependent over the rationals.

In the extreme case where the rank of these coefficients over $\QM$ is
one, the frequencies $(\omega_1,\dots,\omega_n)$ are called
\emph{completely resonant}, and there exist co-prime integers
$p_1,\dots,p_n$ and $\lambda\in\RM$ such that
$\omega_j = p_j \lambda$. In this case, we say that
$(\omega_1,\dots,\omega_n)$ (or the harmonic oscillator
$\sum_{j=1}^{n}\frac{\omega_{j}}{2}(-\h^{2}\frac{{\partial}^2}{{\partial
    x_{j}^2}}+x_{j}^{2})$) exhibits a $p_1:\cdots : p_n$ resonance.

More generally, if there exist an integer $r\in \{1,\dots,n\}$, a
number $\lambda\in\RM$, and co-prime integers $p_{j_1},\dots,p_{j_r}$
such that $\omega_{j_i} = p_{j_i} \lambda$, then we say that
$(\omega_1,\dots,\omega_n)$ exhibits a $p_{j_1}:\cdots : p_{j_r}$
resonance.

In this paper, we are interested in \emph{near resonances}, which we
define as follows.

\begin{defi}

  \label{definition:re}

  Let
  $ \omega_{\varepsilon}= \left(
    \omega_{1,\varepsilon},\dots,\omega_{n,\varepsilon} \right) $ be
  the frequencies of the harmonic oscillator
  $\hat{H}_{2,\varepsilon}$.  One says that $ \omega_{\varepsilon} $
  admits a near resonance of type
  $p_{j_1} \thickapprox \cdots \thickapprox p_{j_r}$ if the map
  $ \varepsilon\rightarrow \omega_{\varepsilon} $ is $ C^{\infty} $
  and $ \omega_{0} $ admits a resonance of type
  $p_{j_1}:\cdots : p_{j_r}$.

  For example, we say that there is a near resonance relation of the
  form $ p \thickapprox q $, where $ p,q \in \NM^* $, if there exist
  $ j_{0},k_{0} \in \lbrace 1,\dots,n \rbrace $ such that
  $ q \omega_{j_{0},0}= p \omega_{k_{0},0} $.
\end{defi}

This concept of near resonance was introduced in molecular
spectroscopy, where such a phenomemon is extremely common among small
molecules, see~\cite{joyeux-birkhoff}.

In order to better understand the non-quadratic terms
$ K_{j,\varepsilon} $ obtained in Theorem~\ref{theo:BGNF} from the
Hamiltonian $ \hat{P}_{\varepsilon}$ in near resonance, we suggest in
this paper to study the typical case of the near Fermi resonance in
dimension $ 2 $, denoted by $ 1\thickapprox 2 $ (\emph{i.e.}
$\omega_{2,0} = 2\omega_{1,0}$), and we will give an explicit
computation of the first non trivial term of the
$\varepsilon$-Birkhoff normal form.  Physically, the near Fermi
resonance is known to affect the spectroscopy of many molecules, among
which we find famous ones such as carbon dioxide $ \textup{CO}_2 $ and
carbon disulfide $ \textup{CS}_2 $, see~\cite{joyeux}. The physics or
chemical literature on Fermi resonance is enormous; interestingly, in
a very recent work, the $1\thickapprox 2$ resonance of the
$ \textup{CO}_2 $ molecule was demonstrated to have a direct effect on
global warming~\cite{wordsworth-al-24}. From the point of view of
classical mechanics, the Fermi resonance and its bifurcations have
been extensively studied, see for
instance~\cite{cushman-dullin-etc-07,hanssmann20} and references
therein.

Consider on $ L^{2}(\RM^2) $ the semiclassical Schrödinger operator
which is transformed according to Proposition~\ref{prop:preparation}
into a perturbation of the harmonic oscillator $ \hat{H}_{2,0} $,
which we denote by $\hat{P}_\varepsilon(\h)$ again:
\begin{equation}
  \hat{P}_{\varepsilon}( \h)=W_{\varepsilon}(0)+\hat{H}_{2,0}
  +\varepsilon\hat{L}_{2}+ {\varepsilon}^2 \hat M_{2,\varepsilon} + R_\varepsilon (x)+\mathcal{O}(\varepsilon^\infty\norm{x}^2),
\end{equation}
where
\begin{align}
  \hat{H}_{2,0}&= \frac{\omega_{1,0}}{2} \left(-\h^2\frac{\partial^2}{\partial x_{1}^2}+x_{1}^2 \right)+\frac{\omega_{2,0}}{2} \left(-\h^2\frac{\partial^2}{\partial x_{2}^2}+x_{2}^2 \right),\\
  \hat{L}_{2}&= \frac{\omega_{1,1}}{2} \left(-\h^2 \frac{\partial^2}{\partial x_{1}^2}+x_{1}^2 \right)+\frac{\omega_{2,1}}{2} \left(-\h^2 \frac{\partial^2} {\partial x_{2}^2}+x_{2}^2 \right),\\
  \hat{M}_{2,\varepsilon}&= \frac{\tilde{\omega}_{1,2}(\varepsilon)}{2} \left(-\h^2 \frac{\partial^2}{\partial x_{1}^2}+x_{1}^2 \right)+\frac{\tilde{\omega}_{2,2}(\varepsilon)}{2} \left(-\h^2 \frac{\partial^2} {\partial x_{2}^2}+x_{2}^2 \right),
\end{align}
and $ R_\varepsilon(x)=\mathcal{O}(x^3) $ at the origin.  The
associated symbol is
\begin{equation}
  P_{\varepsilon} = W_{\varepsilon}(0)+H_{2,0}+\varepsilon L_{2}+L_\varepsilon',
\end{equation}
where for $ z_{j}=x_{j}+i\xi_{j},j=1,2 $,
\begin{align}
  H_{2,0}
  &=\tfrac{1}{2} \left( \omega_{1,0}\vert z_{1} \vert ^2 +\omega_{2,0} \vert z_{2} \vert ^2 \right),
  & L_{2} &=\tfrac{1}{2}\left(\omega_{1,1}\vert z_{1} \vert ^2 +\omega_{2,1} \vert z_{2} \vert ^2 \right),\\
  L_\varepsilon'
  &= {\varepsilon}^2 M_{2,\varepsilon}+R_\varepsilon (x)+\mathcal{O}(\varepsilon^\infty), \text{ and }
  &M_{2,\varepsilon} &=\tfrac{1}{2}\left(\tilde{\omega}_{1,2}(\varepsilon)\vert z_{1} \vert ^2 +\tilde{\omega}_{2,2}(\varepsilon) \vert z_{1} \vert ^2 \right)\,.
\end{align}
We now consider the associated Taylor series (which we denote by the
same symbols). Thus, $ L_\varepsilon'\in\mathcal{O}_3$, and notice
that $ \varepsilon L_{2}\in \mathcal{D}_3 $ and commutes with
$ H_{2,0} $.  By applying Theorem~\ref{theo:BGNF}, there exist
$ A_{1,\varepsilon}\in \mathcal{O}_3 $ and
$ K_{3,\varepsilon} \in \mathcal{D}_3 $ such that
$ K_{3,\varepsilon}\in\mathcal{O}_2(x,\xi,\h)$ and
\begin{equation}
  e^{i \h^{-1} \ad{A_{1,{\varepsilon}}}} \left( H_{2,0}+{\varepsilon}L_{2}
    +L_{\varepsilon}' \right) = H_{2,0}+K_{3, \varepsilon}+ {\mathcal{O}}_{4},
\end{equation}
hence
\begin{equation}
  {\varepsilon}L_{2}+L_{\varepsilon}'+\frac{i}{\h} \left[
    A_{1,\varepsilon}, H_{2,0}+{\varepsilon}L_{2}+L_{\varepsilon}'
  \right] _W - K_{3, \varepsilon} \in \mathcal{O}_4\,.
\end{equation} 
Let us write the expansion
\begin{equation}
  R_{\varepsilon}(x)=R_{0}(x)+\varepsilon R_{1}(x)+{\varepsilon}^2 R_{2}(x)+\dots,
\end{equation}
where $ R_{j}(x)={\mathcal{O}} (x^3) $, and we can write
$R_0(x) = R_{3,0} + \mathcal{O}(x^4)$ with
$ R_{3,0}\in \mathcal{D}_3 $.  Since
$ \varepsilon^2 M_{2,\varepsilon}\in \mathcal{O}_4 $, we have
$ L_\varepsilon'= R_{3,0} + \mathcal{O}(4)$.  We obtain:
\begin{equation}
  \varepsilon L_{2}+R_{3,0}=\frac{i}{\h} \left[H_{2,0},A_{1,\varepsilon } \right] _W
  +K_{3,\varepsilon}  + \mathcal{O}_4\,.
\end{equation}
Since
$ R_{3,0}\in\mathcal{D}_3 = \ker\left( i\h^{-1}\ad{H_{2,0}} \right)
\oplus \textup{Im} \left( i\h^{-1}\ad{H_{2,0}} \right) $, we can write
\begin{equation}
  R_{3,0}=\tilde{R}_{3,0}+\frac{i}{\h}\left[
    H_{2,0},\tilde{\tilde{R}}_{3,0} \right]_W\,,
  \label{equ:R_30}
\end{equation}
where $ \tilde{R}_{3,0}\in \mathcal{D}_3 $ and commutes with
$ H_{2,0} $, and $ \tilde{\tilde{R}}_{3,0}\in \mathcal{D}_3
$. Consequently,
\begin{equation}
  K_{3,\varepsilon}=\varepsilon L_{2}+\tilde{R}_{3,0} .
\end{equation}
Notice that, because $H_{2,0}$ is quadratic, Equation~\eqref{equ:R_30}
is equivalent to
$ R_{3,0}=\tilde{R}_{3,0}+ \{H_{2,0},\tilde{\tilde{R}}_{3,0} \}$;
therefore, $\tilde{R}_{3,0}$ does not depend on $\h$ nor
$\varepsilon$. Hence, in order to make the condition that
$ \tilde{R}_{3,0} $ commutes with $ H_{2,0} $ explicit, we may use the
known results from the classical normal form, see for
instance~\cite{cushman-dullin-etc-07}. Let us recall the arguments. By
Proposition \ref{prop:freq}, we have
\begin{equation}
  \tilde{R}_{3,0}=\sum _{\vert \alpha \vert +\vert \beta \vert =3 \atop \langle \omega_0,\beta-\alpha \rangle =0}c_{\alpha \beta}^{(3)}  z^\alpha \bar{z}^\beta.
\end{equation}
According to the definition \ref{definition:re}, a near Fermi
resonance for
$\omega_\varepsilon =( \omega_{1,\varepsilon}, \omega_{2,\varepsilon}
)$ means an exact Fermi resonance for
$\omega_0=( \omega_{1,0}, \omega_{2,0} )$, that is:
\begin{equation}
  \label{eq:Fermi}
  \omega_{2,0}=2 \omega_{1,0}\,.
\end{equation} 
Therefore,
\begin{equation}
  \langle \omega_0,\beta- \alpha \rangle =0 \Leftrightarrow
  (\alpha_1+2\alpha_2)=(\beta_1 +2\beta_2),\label{equ:resonance12}
\end{equation}
where $ \alpha=(\alpha_1,\alpha_2)\in \NM^2 $ and
$ \beta=(\beta_1,\beta_2)\in \NM^2. $ Now, to obtain
$\tilde{R}_{3,0}$, it is necessary to look for all monomials of order
$ 3 $ that satisfy the Fermi resonance
relation~\eqref{equ:resonance12}. Thus, $\tilde{R}_{3,0}$ is generated
by the monomials $\bar{z}_1^2 z_2,z_1^2\bar{z}_2 $.  Since
$\tilde{R}_{3,0}$ is real,
\begin{align}
  \tilde{R}_{3,0}
  & = \mu \textup{Re}(z_1^2\bar{z}_2) + \nu \textup{Im}(z_1^2\bar{z}_2)\\
  & =\frac{\mu}{2}(\bar{z}_1^2z_2 + z_1^2\bar{z}_2)
    + \frac{\nu}{2i}(\bar{z}_1^2z_2 - z_1^2\bar{z}_2),
    \quad \mu,\nu \in \RM.
\end{align}
\textbf{Determination of the coefficients $\mu$ and $\nu$}: We Taylor
expand $R_{3,0}$ at the origin:
\begin{align}
  R_{3,0}(x_1,x_2)
  & = \frac{1}{3!} \left[
    \frac{{\partial}^3 R_{3,0}(0)}{{\partial} x_1^3}x_1^3
    +3 \frac{{\partial}^3 R_{3,0}(0)}{{\partial} x_1^2 {\partial} x_2}x_1^2 x_2 \right.
  \\
  & \left. + 3 \frac{{\partial}^3 R_{3,0}(0)}{{\partial}
    x_2^2 {\partial} x_1}x_2^2 x_1
    + \frac{{\partial}^3 R_{3,0}(0)}{{\partial} x_2^3}x_2^3 \right].
\end{align}
Recall that, due to~\eqref{equ:admissible}, $ \tilde{R}_{3,0}$ is
obtained by writing $R_{3,0}$ in the $(z^\alpha\bar{z}^\beta)$ basis
and keeping only the resonant coefficients, \emph{i.e.} the
coefficients of $z_1^2\bar{z}_2$ and $\bar{z}_1^2z_2$.  Since
$ x_j=\frac{1}{2} (z_j + \bar{z}_j)$, we remark first that only the
coefficient of $ x_1^2 x_2 $ in $ R_{3,0}(x_1,x_2) $ contributes to
$ \tilde R_{3,0} $ (by homogeneity considerations), and second, that
the expansion of $R_{3,0}$ in the $(z^\alpha\bar{z}^\beta)$ basis has
only \emph{real} coefficients (which would not be the case for a
general Hamiltonian depending also on $\xi$).  More precisely, we have
\begin{equation}
  x_1^2 x_2 =\frac{1}{8} \left(z_1^2 z_2 +z_1^2 \bar z_2+ \bar z_1^2 z_2 + \bar z_1^2 \bar z_2 + 2 \vert z_1 \vert ^2 z_2 + 2 \vert z_1 \vert ^2 \bar z_2   \right)\,,
\end{equation}
which gives
\begin{align}
  \tilde R_{3,0}
  &= \frac{1}{16} \frac{{\partial}^3 R_{3,0}(0)}{{\partial} x_1^2 {\partial} x_2} \left( z_1^2 \bar z_2+ \bar z_1^2 z_2 \right)\\
  &=\frac{1}{8}\frac{{\partial}^3 R_{3,0}(0)}{{\partial} x_1^2 {\partial} x_2}
    \textup{Re} \left( z_1^2 \bar z_2 \right).
\end{align}
So, $\nu=0$, and we obtain the following result,
\begin{theo}
  \label{theo:mu-reso12}
  In dimension $2$, the quantum $\varepsilon$-Birkhoff-Gustavson
  normal form of the operator
  $ H_{2, \varepsilon}+R_{\varepsilon} (x) $ in near Fermi resonance
  $ 1 \approx 2 $, is equal to
  $ H_{2,0}+K_{3,\varepsilon}+\mathcal{O}_4 $ with
  \begin{equation}
    K_{3,\varepsilon}=\frac{\omega_{1,1}}{2} \varepsilon \vert z_1 \vert ^2
    + \frac{\omega_{2,1}}{2} \varepsilon \vert z_2 \vert ^2
    + \mu \textup{Re} \left( z_1^2 \bar z_2 \right),
  \end{equation}
  where
  \begin{equation}
    \label{equ:mu}
    \mu = \frac{1}{8}
    \frac{{\partial}^3 R_{3,0}(0)}{{\partial} x_1^2 {\partial} x_2}.
  \end{equation}
  Moreover, the remainder in $\mathcal{O}_4$ also belongs to
  $\mathcal{O}_2(x,\xi,\h)$.
\end{theo}

In order to obtain a fully satisfactory result for the original
Schrödinger operator~\eqref{equ:schrodinger}, let us now express $\mu$
in terms of the potential $V_0$.

\begin{theo}
  \label{theo:mu_explicit}
  Expressed in terms of the original potential $V_0$, the coefficient
  $\mu$ of Theorem~\ref{theo:mu-reso12} is given as follows.
  \begin{align}
    \mu  =  \frac{1}{8\sqrt{2}\omega_{1,0}^{3/2}(1+z^2)^{\frac32}}
    & \left( z   \partial^3_{x_1^3} V_0(0)
      + (1-2z^2) \partial^3_{x_1^2 x_2}V_0(0) \right. \\
    & ~~ + \left. (z^3-2z) \partial^3_{x_1 x_2^2}V_0(0)  + z^2 \partial^3_{x_2^3}V_0(0) 
      \right)
  \end{align}
  with
  \[
    \omega_{1,0} =
    \frac{1}{\sqrt{3}}\left[\left(\partial^2_{x_1^2}V_0(0) -
        \partial^2_{x_2^2}V_0(0)\right)^2 + 4
      \left(\partial^2_{x_1x_2} V_0(0)\right)^2\right]^{1/4} \,,
  \]
  and
  \[
    z= \frac{2\abs{\partial^2_{x_1x_2}V_0(0)}}{3\omega_{1,0}^2 -
      \partial^2_{x_1^2}V_0(0) + \partial^2_{x_2^2}V_0(0)}\,.
  \]
\end{theo}
\begin{demo}
  According to Proposition~\ref{prop:preparation}, see
  also~\eqref{eq:reste}, we can write
  \[
    (\Lambda_\varepsilon g_\varepsilon \tau_\varepsilon
    V_\varepsilon)(\tilde y) = V_\varepsilon(x_\varepsilon) +
    \sum_{j=1}^n \frac{\omega_{j,\varepsilon}}{2}\tilde y_j^2 +
    R_\varepsilon(\tilde y)\,,
  \]
  where the variable $\tilde y$ is given by the three successive
  change of coordinates:
  \[
    \tilde y = \Lambda_\varepsilon U_\varepsilon^*(x-x_\varepsilon)\,,
  \]
  where we abuse of the notation
  $\Lambda_\varepsilon(y_1,\dots,y_n) = (\sqrt{\omega_{1,\varepsilon}}
  y_1, \dots, \sqrt{\omega_{n,\varepsilon}} y_n)$. Therefore, for all
  $\abs{\alpha}\geq 3$,
  \[
    \partial_{\tilde y}^\alpha R_\varepsilon(\tilde y) =
    \partial_{\tilde y}^\alpha (\Lambda_\varepsilon g_\varepsilon
    \tau_\varepsilon V_\varepsilon)(\tilde y) \,.
  \]
  Since $R_\varepsilon(\tilde y) = R_{3,0}(\tilde y) + \mathcal{O}_4$
  (in the sense of the Taylor series at the origin)
  \[
    \frac{{\partial}^3 R_{3,0}}{{\partial} \tilde{y}_1^2 {\partial}
      \tilde{y}_2}(0) = \frac{{\partial}^3 R_{0}}{{\partial}
      \tilde{y}_1^2 {\partial} \tilde{y}_2}(0) + \mathcal{O}_4=
    \frac{{\partial}^3(\Lambda_0 g_0 \tau_0 V_0)} {{\partial}
      \tilde{y}_1^2 {\partial} \tilde{y}_2}(0) + \mathcal{O}_4\,,
  \]
  and hence the equality holds without the $\mathcal{O}_4$, since the
  quantities do not depend on $x,\xi,\varepsilon,\h$. Moreover,
  $\tau_0=\textup{Id}$, so
  \begin{equation}
    \frac{{\partial}^3 R_{3,0}}{{\partial} \tilde{y}_1^2 {\partial}
      \tilde{y}_2}(0) = \frac{{\partial}^3(\Lambda_0 g_0 V_0)}
    {{\partial} \tilde{y}_1^2 {\partial} \tilde{y}_2}(0) =
    \frac{1}{\omega_{1,0}\sqrt{\omega_{2,0}}} \frac{{\partial}^3(
      g_0 V_0)} {{\partial} {y}_1^2 {\partial} {y}_2}(0)
    \,.\label{equ:R30}
  \end{equation}
  Let us now compute the derivatives of $(g_0V_0)(y) = V_0(U_0 y)$.
  Let $(a,b)$ and $(c,d)$ be an orthonormal basis of eigenvectors of
  the Hessian matrix $V_0''(0)$ for the eigenvalues
  $(\omega_{1,0}^2,\, \omega_{2,0}^2)$, so that we can take
  \[
    U_0 =
    \begin{pmatrix}
      a & c\\ b & d
    \end{pmatrix}\,.
  \]
  We have
  \begin{align}
    \frac{{\partial}^3(
    V_0(U_0 y))} {{\partial} {y}_1^2 {\partial} {y}_2}(0)
    &=
      a^2c \partial^3_{x_1^3} V_0(0) + (2abc + a^2d) \partial^3_{x_1^2 x_2}V_0(0)\\
    & + (b^2c + 2abd) \partial^3_{x_1 x_2^2}V_0(0) + b^2d \partial^3_{x_2^3}V_0(0)\,.\label{equ:V3}
  \end{align}

  It remains to express the coefficients $a,b,c,d$, taking into
  account the fact that the eigenvalues of the Hessian matrix
  $V''_0(0)$ are $\omega_{1,0}^2=\lambda$, $\omega_{2,0}^2=4\lambda$;
  one can check that this implies
  \[
    \lambda := \frac{1}{3}\sqrt{\left(\partial^2_{x_1^2}V_0(0) -
        \partial^2_{x_2^2}V_0(0)\right)^2 + 4
      (\partial^2_{x_1x_2}V_0(0))^2 } \,.
  \]

  In fact, it is elementary to check that for a real symmetric
  positive matrix whose eigenvalues are $(\lambda, 4\lambda)$, the
  discriminant of the characteristic polynomial is equal to
  $9\lambda^2$, and the matrix takes the form, for some $x\in\RM$,
  \[
    \begin{pmatrix}
      \frac{x+5\lambda}{2} & \frac{\sqrt{9\lambda^2 - x^2}}{2} \\[1em]
      \frac{\sqrt{9\lambda^2 - x^2}}{2} & \frac{5\lambda - x}{2}
    \end{pmatrix}\,.
  \]
  Its eigenvectors can be written as
  \[
    (1,-z) \quad \text{ and } \quad (z, 1)\,,
  \]
  with $z=\frac{\sqrt{9\lambda^2 - x^2}}{3\lambda-x}\geq 0$. Upon
  normalization, we obtain
  \begin{align}
    a & = \frac{1}{\sqrt{1+z^2}}  &  b =  \frac{-z}{\sqrt{1+z^2}} \\
    c & = \frac{z}{\sqrt{1+z^2}} & d = \frac{1}{\sqrt{1+z^2}}\,.
  \end{align}
  In our case we have
  \[
    x = \partial^2_{x_1^2}V_0(0) - \partial^2_{x_2^2}V_0(0)\,,
  \]
  therefore
  \begin{align}
    z & = \frac{2 \abs{\partial^2_{x_1x_2}V_0(0)}}{3\lambda - x}\\
      & = \frac{2 \abs{\partial^2_{x_1x_2}V_0(0)}}{\left(\left(\partial^2_{x_1^2}V_0(0) -
        \partial^2_{x_2^2}V_0(0)\right)^2 + 4 \left(\partial^2_{x_1x_2}
        V_0(0)\right)^2\right)^{1/2} - \partial^2_{x_1^2}V_0(0) + \partial^2_{x_2^2}V_0(0)}\,,
  \end{align}
  and
  \begin{align}
    a^2c & = \frac{z}{(1+z^2)^{\frac32}} \\
    2abc+ a^2d & =\frac{1-2z^2}{(1+z^2)^{\frac32}}\\
    b^2c + 2 abd &= \frac{z^3 - 2 z}{(1+z^2)^{\frac32}}\\
    b^2d &= \frac{z^2}{(1+z^2)^{\frac32}}\,,
  \end{align}
  which, in view of~\eqref{equ:mu},~\eqref{equ:R30}, and
  \eqref{equ:V3}, proves the theorem.
\end{demo}

\section{Spectral analysis in near Fermi resonance}
\label{sec:spectr-analys-fermi}

Theorem~\ref{theo:BGNF} gives a formal conjugation of the initial
Schrödinger Hamiltonian $\hat P_\varepsilon$ into an operator of the
form $c_\varepsilon+ \hat H_{2,0} + K_\varepsilon $, where
$c_\varepsilon=W_\varepsilon{(0)} = \min V_\varepsilon$, see
Proposition~\ref{prop:preparation}. Therefore, it is expected that, in
regimes where $E$ and $\h$ are small enough, the spectrum inside
$(-\infty,E)$ of the normal form
$c_\varepsilon + \hat H_{2,0} + K_\varepsilon $, restricted to the
spectral subspace of $H_{2,0}$ corresponding to energies in
$[0,(1+\eta)E]$, where $\eta>0$ is fixed, is a good approximation of
the spectrum of $\hat P_\varepsilon$ in the interval $(-\infty,E)$,
see~\cite{san-charles}.

The goal of this section is to describe the spectrum of the
$\varepsilon$-Birkhoff-Gustavson normal form of $\hat P_\varepsilon$
in the case of a near Fermi resonance, in terms of the original
potential $V_\varepsilon$. We shall use the expression of the normal
form modulo $\mathcal{O}_4$ given in the previous section. Since the
remainder also belongs to $\mathcal{O}_2(x,\xi,\h)$ by
Theorem~\ref{theo:mu-reso12}, it belongs to
\[
  \mathcal{O}_4(x,\xi,\h) + \varepsilon\mathcal{O}_3(x,\xi,\h) +
  \varepsilon^2\mathcal{O}_2(x,\xi,\h)\,.
\]
Therefore, using~\cite[Lemma 4.2]{san-charles}, the normal form at
order 3, that is to say
$c_\varepsilon+ \hat H_{2,0} + K_{3,\varepsilon} $, is expected to
approximate the spectrum of $\hat P_\varepsilon$ below some energy $E$
with a precision of order
$\mathcal{O}(E^2)+ \mathcal{O}(\varepsilon E^{3/2}) +
\mathcal{O}(\varepsilon^2 E)$. For instance, if we are interested in a
fixed number of low-lying eigenvalues, we may take $E=C\h$ for some
constant $C>0$ and we obtain a precision of order
$\mathcal{O}(\h^2)+ \mathcal{O}(\varepsilon \h^{3/2}) +
\mathcal{O}(\varepsilon^2 \h)$.

In order to compute the matrix elements of the normal form, we shall
need to understand the Weyl quantization $ \hat{K}_{3,\varepsilon} $
of $ K_{3,\varepsilon} $, which is
\begin{align}
  \hat{K}_{3,\varepsilon}
  & =\varepsilon \hat{L}_{2}+\mu \widehat{\textup{Re}(z_1^2 \bar z_2 )},\\
  & = \varepsilon \hat{L}_{2}+\mu
    \left( x_1^2 x_2 -2 \h ^2 x_1 \frac{\partial}{\partial x_1} \frac{\partial}{\partial x_2}
    + \h ^2 x_2 \frac{\partial ^2}{\partial x_1^2}
    - \h^2\frac{\partial}{\partial x_2 } \right).
\end{align} 
For that purpose, it will be convenient to pass to the Bargmann
representation.

\subsection{Creation and annihilation operators}
Let us introduce
\begin{equation}
  \begin{array}{c}
    a_j(\h)= \dfrac{1}{\sqrt{2\h}} \left(x_j+\h \dfrac{\partial}{\partial x_j} \right); \\
    b_j(\h)= \dfrac{1}{\sqrt{2\h}} \left(x_j-\h \dfrac{\partial}{\partial x_j} \right),
  \end{array}\label{equ:ab}
\end{equation}
which are respectively called the operators of annihilation and
creation, acting as unbounded operators on $L^2(\mathbb{R}^n)$ (see
for instance~\cite{combescure-robert-book}).
 
Formally, the operators $a_j(h)$ and $ b_j(h)$ satisfy the following
properties:
\begin{equation}
  \begin{array} {cll}
    a_j^*(\h)= b_j(\h), & b_j^*(\h)=  a_j(\h),  \\
    \left[ a_j(\h), b_k(\h) \right] = \delta_{jk}, & \left[ a_j(\h), a_k(\h) \right]=0, & \left[ b_j(\h), b_k(\h) \right]=0. \\ 
  \end{array}
\end{equation}
Using $a_j(\h)$ and $b_j(\h)$ to rewrite $\hat{H}_{2,0}$, we have
\begin{equation}
  \hat{H}_{2,0} = \h\left( \omega_{1,0}(a_1(\h) b_1(\h)-\tfrac{1}{2})
    + \omega_{2,0}(a_2(\h) b_2(\h)-\tfrac{1}{2}) \right),
\end{equation}
\begin{equation}
  \hat{L}_{2} =
  \h\left( \omega_{1,1}(a_1(\h) b_1(\h)-\tfrac{1}{2})
    + \omega_{2,1}(a_2(\h) b_2(\h)-\tfrac{1}{2}) \right),
\end{equation}
and $\hat{K}_{3,\varepsilon}$ is given by
\begin{equation} \label{creator K3}
  \begin{array}{cll}
    \hat{K}_{3,\varepsilon}
    &= &\h\varepsilon \left( \omega_{1,1}(a_1(\h) b_1(\h)-\tfrac{1}{2})
         +\omega_{2,1}(a_2(\h) b_2(\h)-\tfrac{1}{2}) \right)\\
    & & \\
    & & + \sqrt{2} \mu \h^{3/2} \left( a_2(\h)b_1^2(\h)+a_1^2(\h)b_2(\h) \right).
  \end{array} 
\end{equation}

\subsection{Bargmann representation}
\label{sec:bargm-repr}
 
In this section, we review several fundamental results relating to the
Bargmann transform and the space $ \mathcal{B_F} $ of Bargmann-Fock,
also known as the Bargmann space~\cite{bargmann}. For further details
see~\cite{combescure-robert-book}. Let us consider the space
\begin{equation}
  \mathcal{B_F}=\left\{ \varphi(\zeta) \text{ holomorphic function on} \ \mathbb{C}^n ; \int_{\mathbb{C}^n} |\varphi(\zeta) |^2 d\mu_n(\zeta) < + \infty  \right\},
\end{equation}
where $ \dd \mu_n(\zeta)$ is the Gaussian measure defined by
$ \dd \mu_n(\zeta) = e^{-|\zeta|^2/\h} \dd{}^n \zeta $,
$\zeta_j = \frac{x_j - i\xi_j}{\sqrt{2}} = \bar{z_j}/\sqrt{2}$, and
$\dd{}^n \zeta$ is the Lebesgue measure in $(x,\xi)$.
$\mathcal{B_F} $ is a Hilbert space equipped with natural inner
product
\[
  \pscal{f}{g} = \int_{\mathbb{C}^n}
  f(\zeta)\overline{g(\zeta)}d\mu_n(\zeta).
\]

The semiclassical scaling that we use here can be found, for instance,
in~\cite[Proposition 7]{combescure-robert-book}.
\begin{theo}[\cite{bargmann}]
  There is a unitary mapping $\mathrm{T_\mathcal{B}}$ from
  $ L^{2}(\mathbb{R}^n) $ to $\mathcal{B_F} $ defined by
  \begin{equation}
    (\mathrm{T_\mathcal{B}} f)(\zeta)= \dfrac{1}{2^{\frac{n}{2}}(\pi \h )^{\frac{3n}{4}}}
    \int_{\mathbb{R}^n} f(x)e^{-\frac{1}{\h}[\frac{1}{2}(\zeta^2+x^2)+\sqrt{2}x\zeta]} \dd{x} , 
  \end{equation}
  $ \mathrm{T_\mathcal{B}} $ is known as the Bargmann transform.
\end{theo}

The Bargmann representation is particularly suited for studying
harmonic approximations, as creation and annihilation operators become
expressible in a very simple way.
\begin{proposition}[\cite{bargmann}]
  \label{prop:ab}
  We have
  \begin{equation}
    \mathrm{T_\mathcal{B}}(a_j(\h))\mathrm{T_\mathcal{B}}^{-1}=
    \deriv{}{\zeta_j} \quad \text{and} \quad
    \mathrm{T_\mathcal{B}}(b_j(\h))\mathrm{T_\mathcal{B}}^{-1}= \zeta_j,
  \end{equation}
  where $\zeta_j$ represents the operator of multiplication by
  $\zeta_j$.
\end{proposition}
Thus, the harmonic oscillator's Bargmann transform is,
\begin{equation}
  \hat{H}_{2,0}^\mathcal{B} :=
  \mathrm{T_\mathcal{B}}(\hat{H}_{2,0})\mathrm{T_\mathcal{B}}^{-1}= \h \sum_{j=1}^{n}\omega_{j,0}
  \left( \zeta_j \frac{\partial}{\partial \zeta_j} + \frac{1}{2}\right).
\end{equation}
It is then standard to prove that the functions
$\left\{ \phy_\alpha := \frac{\zeta^\alpha}{
    (2\pi\h)^{n/2}\sqrt{\alpha !}} \right\}_{ \alpha \in \mathbb{N}^n}
$ form an orthonormal Hilbertian basis of $ \mathcal{B_F}$, on which
$\hat{H}_{2,0}^{\mathcal{B}}$ is diagonal. It follows that
$\hat{H}_{2,0}^{\mathcal{B}}$ is essentially self-adjoint with a
discrete spectrum which consists of the eigenvalues
$ \lambda_\alpha $, $\alpha\in\NM^n$ given by:
\[
  \sigma ( \hat{H}_{2,0}^{\mathcal{B}}) = \Big\{ \lambda_\alpha = \h
  \big(u+ \frac{\abs{\omega_0}}{2}\big), \ \ u =
  \pscal{\omega_0}{\alpha} = \sum_{j=1}^{n} \omega_{j,0} \alpha_j, \,
  \alpha\in\NM^n\Big\},
\]
(repeated with their multiplicities) where
$\omega_0 = (\omega_{1,0},...,\omega_{n,0})$,
$\alpha = (\alpha_1,...,\alpha_n)$ and
$ \abs{\omega_0} = \sum_{j=1}^{n} \omega_{j,0}$.  The corresponding
eigenspaces are
\[
  \mathcal{H}_u^{\mathcal{B}} = \textup{Span} \left\{
    \frac{\zeta^\alpha}{ (2\pi\h)^{n/2} \sqrt{\alpha !}}, \,
    \pscal{\omega_0}{\alpha}=u \right\}.
\]

\subsection{Spectrum in near Fermi resonance }
    
The authors in~\cite{san-charles} give a detailed study of the
spectrum of semiclassical pseudodifferential operators whose harmonic
approximation possesses an exact complete resonance; the idea was to
restrict the Birkhoff normal form to the eigenspaces of the resonant
harmonic oscillator. In this section, we extend this result to the
case of a near Fermi resonance, using the
$\varepsilon$-Birkhoff-Gustavson normal form
$H_{2,0} + K_{3,\varepsilon} + \mathcal{O}(4)$ of
Theorem~\ref{theo:BGNF}.

First, let us rewrite this normal form in the Bargmann
representation. According to (\ref{creator K3}), we have
\begin{equation}\label{equ:barg_K3}
  \begin{array}{cll} 
    \hat{K}_{3,\varepsilon}^{\mathcal{B}} & := & \mathrm{T_\mathcal{B}}(\hat{K}_{3,\varepsilon})\mathrm{T_\mathcal{B}}^{-1}=\h\varepsilon \omega_{1,1}
                                                 \left( \zeta_1 \frac{\partial}{\partial \zeta_1}+\tfrac{1}{2}\right)+\h\varepsilon \omega_{2,1}\left( \zeta_2 \frac{\partial}{\partial \zeta_2}+\tfrac{1}{2}\right) \\
                                          && +\sqrt{2} \mu \h^{3/2} \left(\zeta_2\frac{\partial^2}{\partial \zeta_1^2}+\zeta_1^2 \frac{\partial}{\partial \zeta_2} \right). 
  \end{array}
\end{equation}

Thus, in order to obtain a good approximation of the spectrum of our
Schrödinger operator $ \hat{P}_\varepsilon $, we shall compute the
spectrum of the restriction of
$ \hat{K}_{3,\varepsilon}^{\mathcal{B}}$ to the eigenspaces of
$ \hat{H}_{2,0}^{\mathcal{B}}$ (recall that
$ \hat{H}_{2,0}^{\mathcal{B}}$ commutes with
$ \hat{K}_{3,\varepsilon}^{\mathcal{B}} $). To this aim, we calculate
the matrix elements of $\hat{K}_{3,\varepsilon}^{\mathcal{B}}$ in the
basis $\phy_\alpha$ of $\mathcal{H}_N^\mathcal{B}$, where $N\in\NM$
and
\begin{equation}
  \label{equ:eigenspace}
  \mathcal{H}_N^\mathcal{B}= \textup{Span} \left\{\varphi_\alpha=
    \frac{\zeta^{\alpha_1}\zeta^{\alpha_2} }{2\pi\h \sqrt{\alpha_1 !\alpha_2 !}} ; \quad   \alpha=(\alpha_1 , \alpha_2) \in \NM^2, \,  \alpha_1 + 2 \alpha_2 = N \right\} .
\end{equation} 

As in~\cite{san-kaoutar13}, we use
\begin{align}\label{equ:action_dz}
  \frac{\partial \varphi_{(\alpha_1,\alpha_2)}}{\partial \zeta_1}
  &=\sqrt{\alpha_1}\varphi_{(\alpha_1-1,\alpha_2)}\,, \\ 
  \frac{\partial^2 \varphi_{(\alpha_1,\alpha_2)}}{\partial \zeta_1^2}
  &= \sqrt{\alpha_1}\sqrt{\alpha_1-1}\varphi_{(\alpha_1-2,\alpha_2)}\,, \\
  \frac{\partial \varphi_{(\alpha_1,\alpha_2)}}{\partial \zeta_2}
  &= \sqrt{\alpha_2}\varphi_{(\alpha_1,\alpha_2-1)}, \\
  \frac{\partial^2 \varphi_{(\alpha_1,\alpha_2)}}{\partial \zeta_2^2}
  &= \sqrt{\alpha_2}\sqrt{\alpha_2-1}\varphi_{(\alpha_1,\alpha_2-2)}\,,\\
  \frac{\partial^2 \varphi_{(\alpha_1,\alpha_2)}}{\partial \zeta_1 \partial \zeta_2}
  &=\sqrt{\alpha_1} \sqrt{\alpha_2}\varphi_{(\alpha_1-1,\alpha_2-1)}\,,
\end{align}
and
\begin{align}
  \zeta_1\varphi_{(\alpha_1,\alpha_2)}
  &= \sqrt{\alpha_1+1}\varphi_{(\alpha_1+1,\alpha_2)}, \\ 
  \zeta_1^2\varphi_{(\alpha_1,\alpha_2)}
  &= \sqrt{\alpha_1+2}\sqrt{\alpha_1+1}\varphi_{(\alpha_1+2,\alpha_2)},\\ 
  \zeta_2\varphi_{(\alpha_1,\alpha_2)}
  &=  \sqrt{\alpha_2+1}\varphi_{(\alpha_1,\alpha_2+1)},\\ 
  \zeta_2^2\varphi_{(\alpha_1,\alpha_2)}
  &= \sqrt{\alpha_2+2}\sqrt{\alpha_2+1}\varphi_{(\alpha_1,\alpha_2+2)}, \\ 
  \zeta_1\zeta_2\varphi_{(\alpha_1,\alpha_2)}
  &= \sqrt{\alpha_1+1}\sqrt{\alpha_2+1}\varphi_{(\alpha_1+1,\alpha_2+1)}.
\end{align}
According to~\eqref{equ:barg_K3}, we obtain, with
$\alpha=(\alpha_1,\alpha_2)$,
\begin{align}
  \hat{K}_{3,\varepsilon}^{\mathcal{B}}\varphi_{\alpha}
  =\ & 
       \h\varepsilon \omega_{1,1} \left(\zeta_1 \partial_{\zeta_1}
       \phy_\alpha +\tfrac{1}{2} \phy_\alpha \right) + \h \varepsilon \omega_{2,1} \left(\zeta_2 \partial_{\zeta_2} \phy_\alpha
       + \tfrac{1}{2} \phy_\alpha \right) \\ 
     & +\sqrt{2}\mu \h^{3/2} \left(\zeta_2\sqrt{\alpha_1 }\sqrt{\alpha_1-1}\varphi_{(\alpha_1-2,\alpha_2)}+\zeta_1^2\sqrt{\alpha_2}\varphi_{(\alpha_1,\alpha_2-1)}\right), \\ 
  =\ & \h \left(\varepsilon \omega_{1,1} \left( \alpha_1+\tfrac{1}{2}\right)+\varepsilon \omega_{2,1} \left( \alpha_2+\tfrac{1}{2}\right)\right)\varphi_{\alpha} \\
     &   +\sqrt{2}\mu \h^{3/2} \left(\sqrt{(\alpha_2+1)\alpha_1(\alpha_1-1)}\varphi_{(\alpha_1-2,\alpha_2+1)} \right.\\
     &   \left. + \sqrt{(\alpha_1+2)(\alpha_1+1)\alpha_2} \varphi_{(\alpha_1+2,\alpha_2-1)} \right).
\end{align}
We can confirm on these formulas that the space
$\mathcal{H}_N^\mathcal{B}$ is stable by
$\hat{K}_{3,\varepsilon}^{\mathcal{B}}$ because:
\begin{equation}
  \alpha_1-2+2(\alpha_2+1)=\alpha_1+2\alpha_2=N \ \text{and} \ \alpha_1+2+2(\alpha_2-1)=\alpha_1+2\alpha_2=N\,.
\end{equation}
We can also confirm the Hermitian symmetry of the matrix
$\hat{K}_{3,\varepsilon}^{\mathcal{B}}$. Indeed, let us order the
basis $\phy_\alpha$ according to
$\ell=0,1,...,E \left[ \frac{N}{2}\right]$ and
$\alpha_\ell=(N-2\ell, \ell)$. The image of $\phy_{\alpha_\ell}$ is a
vector with three components, on $\phy_{\alpha_\ell}$,
$\phy_{\alpha_{\ell-1}}$, and $\phy_{\alpha_{\ell+1}}$. Thus, the
matrix is symmetric if and only if the coefficient of
$\hat{K}_{3,\varepsilon}^{\mathcal{B}}(\phy_{\alpha_\ell})$ on
$\phy_{\alpha_{\ell+1}}$ is equal to the coefficient of
$\hat{K}_{3,\varepsilon}^{\mathcal{B}}(\phy_{\alpha_{\ell+1}})$ on
$\phy_{\alpha_\ell}$. The first one is equal to
\[
  \sqrt{2}\mu \h^{3/2} \sqrt{(\alpha_2+1)\alpha_1(\alpha_1-1)}
\]
with $\alpha_1=\ell$, $\alpha_2=N-2\ell$, and the second one is
\[
  \sqrt{2}\mu \h^{3/2} \sqrt{(\alpha'_1+2)(\alpha'_1+1)\alpha'_2}
\]
with $\alpha'_1 = \alpha_1-2$, $\alpha'_2=\alpha_2+1$, so the equality
of these coefficients is clear.

In matrix form, we obtain
\begin{theo}
  \label{theo:matrix}
  For any $N\geq 0$, in the basis $(\phy_\alpha)$ of the space
  $ \mathcal{H}_N^\mathcal{B}$ given in~\eqref{equ:eigenspace},
  ordered by $\alpha_2$, of size $\lfloor \frac{N}{2} \rfloor + 1$,
  the matrix of the operator $\hat{K}_{3,\varepsilon}^{\mathcal{B}}$
  is
  \begin{equation}
    \mathcal{M}_{\varepsilon, N} =
    \begin{pmatrix}
      d_{N, 0} & A_{N, 0} & & & &  & & \\
      A_{N, 0} & d_{N, 1} & \ddots & & & \mathbf{0} & \\
      & \ddots & \ddots & \ddots &  &   \\
      & & A_{N, \ell-1} & d_{N, \ell} & A_{N, \ell} &  & \\
      & \mathbf{0} & & A_{N, \ell} & d_{N, \ell+1} & \ddots &  \\
      & & & & A_{N, \ell+1} & \ddots & \ddots\\
      & & & &  & \ddots & \ddots\\

   \end{pmatrix}
   ,
 \end{equation}
 where for $\ell =0,1,..., \lfloor\frac{N}{2}\rfloor$:
 \begin{equation}
   \left\{
     \begin{aligned}
       A_{N, \ell}=& \sqrt{2} \mu \h^{3 / 2} \sqrt{(\ell+1)(N-2 \ell)(N-2\ell-1)} \\
       d_{N,\ell}=& \h\varepsilon \omega_{1,1} \left(
         N-2\ell+\tfrac{1}{2}\right)+ \h \varepsilon \omega_{2,1}
       \left( \ell+\tfrac{1}{2}\right),
     \end{aligned}
   \right.
 \end{equation}
 and $\mu$ is defined in Theorem~\ref{theo:mu_explicit}.
\end{theo}

Finally, note that the spectrum of
$ \hat{H}_{2,0}^{\mathcal{B}} + \hat{K}_{3,\varepsilon}^{\mathcal{B}}$
restricted to $ \mathcal{H}_N^\mathcal{B}$ is exactly
\[
  \bigcup_{N\in\NM}\left\{\h\omega_{1,0}(N+ \tfrac{3}{2}) +
    \mu_{\varepsilon,N,j}; \quad j=0,\dots \lfloor \tfrac{N}{2}
    \rfloor\right\}\,,
\]
where
$(\mu_{\varepsilon,N,0},\dots,\mu_{\varepsilon,N, \lfloor \frac{N}{2}
  \rfloor})$ are the eigenvalues of the matrix
$\mathcal{M}_{\varepsilon, N}$ of Theorem~\ref{theo:matrix}.

\section{Numerical simulations for the near Fermi resonance}
\label{sec:numerics}

We now illustrate our results numerically for the following family of
potentials:
\[
  V_\varepsilon(x_1,x_2) = (\tfrac{1}{2}+\varepsilon)x_1^2 +
  (2+2c\varepsilon) x_2^2+ x_1^3 + \tfrac{1}{2}x_1 x_2^2 + \gamma
  x_1^2 x_2 + x_1^4 + x_2^4\,,
\]
where $c\in\RM$ and $\gamma$ are fixed, and $\varepsilon\in\RM$ is
such that
\begin{equation}
  \tfrac{1}{2}+\varepsilon >0 \quad \text{ and } 1+c\varepsilon >0 \,.
  \label{equ:V-min-non-deg}
\end{equation}
In view of the normalization~\eqref{equ:Wepsilon}, this potential
corresponds to $W_\varepsilon(0)=0$,
$\omega_{1,\varepsilon}:=\sqrt{1+2\varepsilon}$ and
$\omega_{2,\varepsilon}=2\sqrt{1+c\epsilon}$. Therefore,
\[
  \omega_{1,0}=1, \omega_{2,0}=2, \omega_{1,1}=1, \omega_{2,1}=c\,.
\]
Since the quadratic part of the potential is already diagonalized, we
can take $U=\textup{Id}$, and hence it is easy to compute $z=0$ and
\begin{equation}
  \mu=\frac{\gamma}{4\sqrt{2}}\label{equ:mu-exemple}
\end{equation}
from Theorem~\ref{theo:mu_explicit}.

Due to the term $ x_1^4 + x_2^4$, the potential $V_\varepsilon$ is
confining, and hence
$\hat P_\varepsilon:=-\frac{\h^2}{2}\Delta+V_\varepsilon$ has a
discrete spectrum, bounded from below, and whose eigenvalues form a
non-decreasing sequence tending to $+\infty$. We are interested in the
eigenvalues that belong to the interval $(-\infty, E)$, where $E\to 0$
as $\h\to 0$.

\subsection{Numerical computation of the spectrum of
  $\hat P_\varepsilon$}
\label{sec:numer-comp-spectr}

In order to numerically compute the spectrum of $\hat P_\varepsilon$
without any \emph{a priori}, we use a general spectral method, not
specially adapted to the $1:2$ resonance. Namely, we consider the
Hermite functions $H_{\alpha_1,\alpha_2}(x_1,x_2)$, which correspond
via the Bargmann transform to the functions $\phy_{\alpha_1,\alpha_2}$
from~\eqref{equ:eigenspace}. We order them according, first, to the
energy level $\h(\alpha_1+\alpha_2+1)$ of the usual harmonic
oscillator $-\frac{\h^2}{2}\Delta + \frac{1}{2}(x_1^2 + x_2^2)$, and
then according to $\alpha_1$. In other words, the first few pairs
$(\alpha_1,\alpha_2)$ in increasing order are:
\[
  (0,0), (0,1), (1,0), (0,2), (1,1), (2,0), (0,3), (1,2), (2,1),
  (3,0), \text{ etc.}
\]
The next step of the numerical method is to truncate the basis
$(H_\alpha)$ into a finite family adapted to the chosen spectral
window; according to the semiclassical theory, for a given $E$,
eigenfunctions corresponding to eigenvalues less than $E$ must be
microlocalized in the phase space region $\Omega_\varepsilon$ given by
$\frac{1}{2}\norm{\xi}^2+V_\varepsilon(x) \leq E$. Assume that the
parameters $\gamma,c$ are chosen such that $V_\varepsilon$ admits a
unique minimum at $x=0$, for $\varepsilon$ small enough (see
Proposition~\ref{prop:unique-min} below). By construction, we have
\[
  \frac{1}{2}\norm{\xi}^2+V_\varepsilon(x) =
  \frac{1}{2}(x_1^2+\xi_1^2) + 2x_2^2+\frac{1}{2}\xi_2^2 +
  \mathcal{O}(\varepsilon x^2) + \mathcal{O}(x^3).
\]
Therefore, if $E$ and $\varepsilon$ are small enough,
$\Omega_\varepsilon$ is contained in the ball of radius $\sqrt{3E}$ in
phase space. (The number $3$ could be replaced by any number larger
than $2$. For large $E$, one could do better by using the confinement
gained by the quartic term $x^4$, but in this work we are mainly
interested in small $E$). Taking into account that a wave function
localized inside a ball may extend slightly beyond the ball, typically
at a distance of order $\mathcal{O}(\sqrt{\h})$, we can decide to
truncate the basis for $(\alpha_1,\alpha_2)$ satisfying
$\h(\alpha_1+\alpha_2+1) \leq \frac{3}{2}E + \sqrt{\h}$; thus we shall
take $\alpha_1+\alpha_2 \leq M$ with
\begin{equation}
  M = \frac{3E}{2\h}+\frac{1}{\sqrt{\h}}\,.\label{equ:M}
\end{equation}
We obtain a basis $\mathcal{B}_M$ of cardinal $(M+1)(M+2)/2$ which,
when $E$ is bounded and $\h$ is small, is
$\mathcal{O}(\frac{E^2}{\h^2})$.

On this basis, the action of differential operators with polynomial
coefficients can be explicitly computed, similarly
to~\eqref{equ:action_dz}; namely, from~\eqref{equ:ab} and
Proposition~\ref{prop:ab}, we have
\[
  x_j H_\alpha = \sqrt{\frac{\h}{2}}\left(
    \sqrt{\alpha_j}H_{\alpha-1_j} + \sqrt{\alpha_j+1}
    H_{\alpha+1_j}\right)
\]
and
\[
  \h\partial_{x_j} H_\alpha = \sqrt{\frac{\h}{2}}\left(
    \sqrt{\alpha_j}H_{\alpha-1_j} - \sqrt{\alpha_j+1}
    H_{\alpha+1_j}\right)\,,
\]
where we have denoted by $(1_1,1_2)$ the canonical basis of $\ZM^2$.
In order to implement the operator of multiplication by $x_j^2$,
instead of writing the explicit formula, we may also simply square the
truncated matrix for $x_j$; by doing so, of course we introduce an
error due to the fact that matrix truncation does not commute with
matrix multiplication. The wrong columns concern the images of
$H_\alpha$ when $\alpha+1_j\not\in[0,M] \times [0,M]$, \emph{i.e.}
$\alpha_1+\alpha_2+1 > M$. In other words, in order to obtain a
correct matrix, we need to delete the ``highest energy level'', which
means truncating the computed matrix to the smaller basis
$\mathcal{B}_{M-1}$. Similarly, since the matrix for $x_j^2$ has a
band structure of width 5 (instead of 3 for $x_j$), when computing
$x_j^3$ we need to reduce $M$ by 2 more. Finally, the truncated matrix
for $x_j^4$ will be exact if we reduce $M$ by 3 more. The analogous
discussion holds for $\h\partial_{x_j}$. Consequently, in order to
obtain the matrix for $\hat{P}_\varepsilon$ on $\mathcal{B}_M$, we
need to start from the larger basis $\mathcal{B}_{M+6}$.

In view of the above discussion, we can now implement, explicitly, the
matrix of $\hat{P}_\varepsilon$ on the basis
$\mathcal{B}_M=\{H_\alpha, \abs{\alpha}\leq M\}$, and call a standard
diagonalization routine for real symmetric matrices. On a standard
laptop, this can be easily done for a matrix of size
$1000 \times 1000$, which means $M\leq 43$.

\begin{figure}[h]
  \centering
  \begin{subfigure}[b]{0.45\linewidth}
    \includegraphics[width=\linewidth]{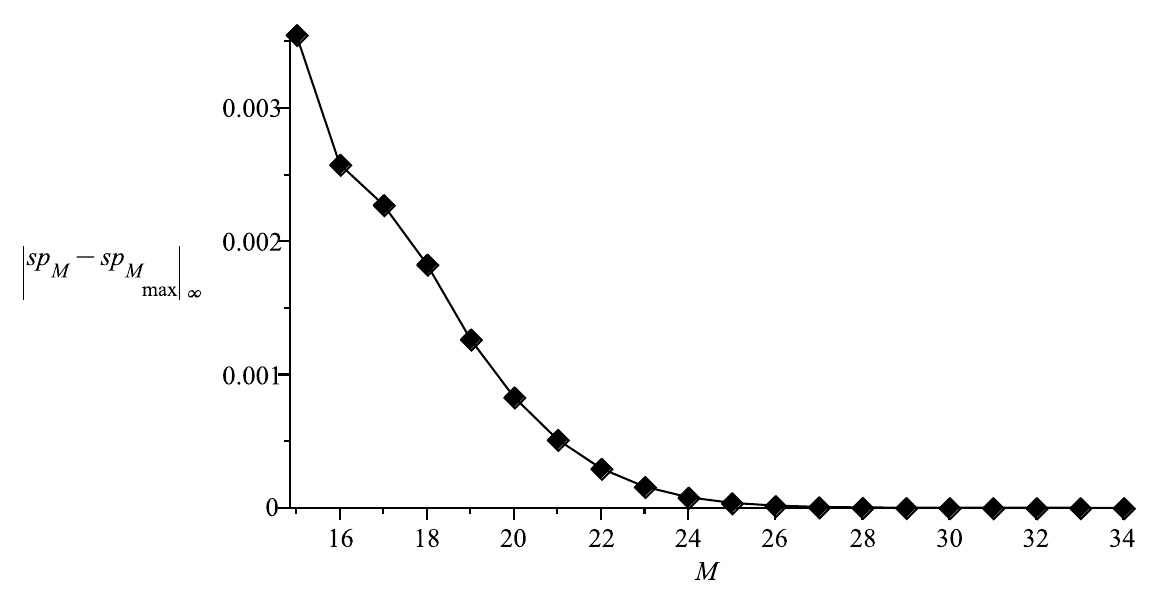}
    \caption{$\norm{\textup{sp}_M -
        \textup{sp}_{M_{\textup{max}}}}_{\infty}$}
  \end{subfigure}
  \begin{subfigure}[b]{0.45\linewidth}
    \includegraphics[width=\linewidth]{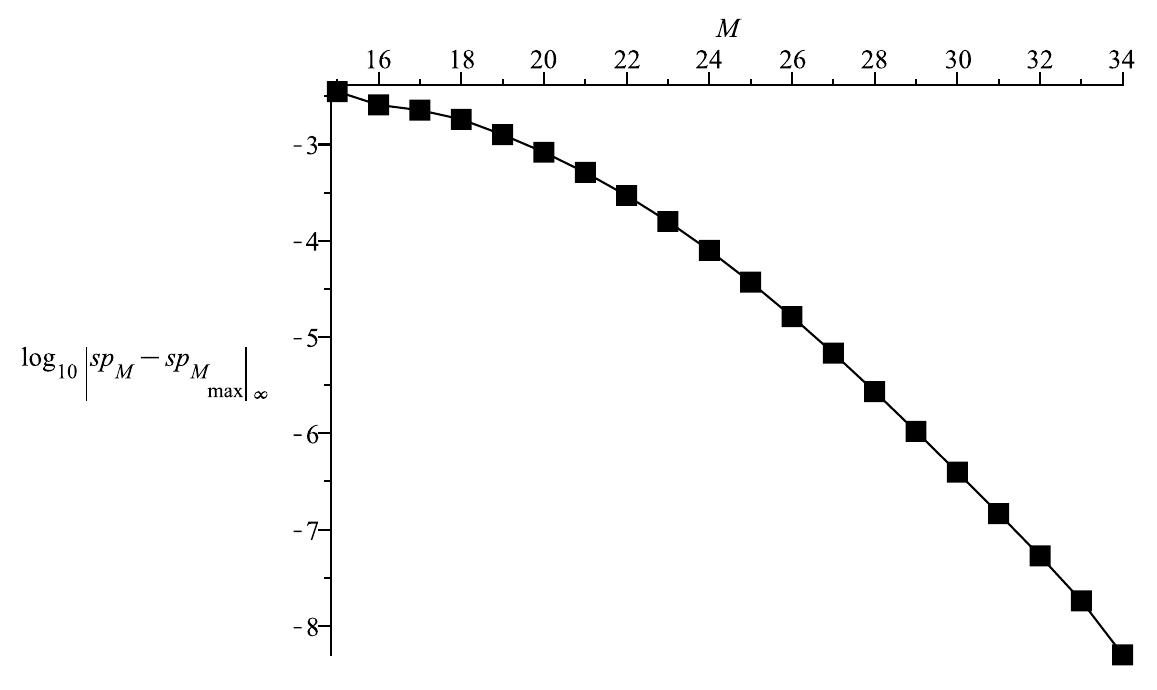}
    \caption{$\log_{10}(\norm{\textup{sp}_M -
        \textup{sp}_{M_{\textup{max}}}}_{\infty})$}
  \end{subfigure}
  \caption{Quality of the numerical spectrum in terms of $M$. With
    $\h=0.01$, $\varepsilon=0$, and $E=10\h$, the value given
    by~\eqref{equ:M} is $M=25$. We plot here
    $\norm{\textup{sp}_M - \textup{sp}_{M_{\textup{max}}}}_{\infty}$
    with $M_{\textup{max}} = 35$, and $\textup{sp}_M$ is the spectrum
    below $E$ obtained by diagonalizing the truncated matrix of
    $\hat{P}_\varepsilon$ in the basis $\mathcal{B}_M$. The same data
    with a $\log_{10}$ scale is plotted on the right picture. We see
    that, at $M\geq 25$, the error is indeed
    negligible. Interestingly, we also see that simply choosing
    $M=\frac{3E}{2\h}=15$ would \emph{not} be sufficient.}
  \label{fig:erreurM}
\end{figure}

Our first experiment is to test the validity of the
truncation~\eqref{equ:M}. Taking $\h=0.01$ and $E=10\h$,
Formula~\eqref{equ:M} gives $M=25$, which corresponds to a matrix of
size $351\times 351$. In Figure~\ref{fig:erreurM}, we vary $M$ around
that value, and plot the $\ell_\infty$ norm of the difference between
the spectrum below $E$ computed with the given $M$ value (that is,
using the matrix obtained from the basis $\mathcal{B}_M$) and the
spectrum computed using the largest $M$ (here
$M_{\textup{max}}=35$). This experiment confirms that the value
predicted by~\eqref{equ:M} is large enough to obtain a very good
accuracy.

\begin{figure}[h]
  \centering \includegraphics[width=0.7\linewidth]{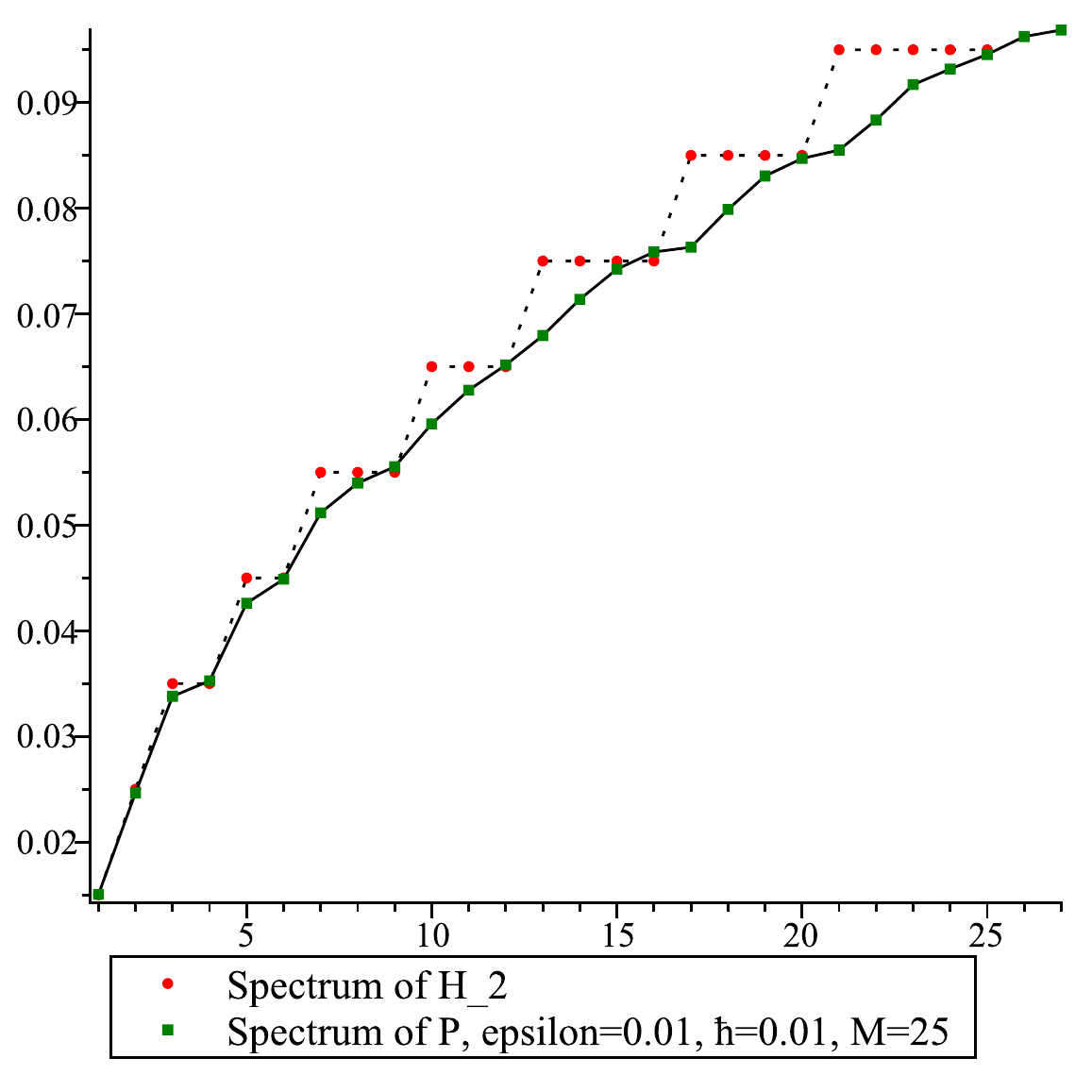}
  \caption{Spectrum of $\hat{P}_\varepsilon$ (green boxes, solid line)
    on top of the spectrum of the $1:2$ harmonic oscillator $H_{2,0}$
    (red discs, dotted line).}
  \label{fig:spectrum_P_H2}
\end{figure}

In Figure~\ref{fig:spectrum_P_H2}, we compare the spectrum of
$\hat{P}_\varepsilon$ obtained with the above numerical method to the
spectrum of the unperturbed harmonic oscillator $H_{2,0}$. Recall that
the oscillator spectrum is explicit
\begin{equation}
  \textup{sp}(H_{2,0}) =   \left\{\h\omega_{1,0}(\alpha_1 + 2\alpha_2 +
    \tfrac{3}{2}),\quad (\alpha_1,\alpha_2)\in \NM^2\right\}\,
  \label{equ:polyads}
\end{equation}
and features the famous phenomemon of \emph{clustering} of eigenvalues
on the ladder $\h\omega_{1,0}(N+\frac{3}{2}), N\geq 0$, which
correspond to \emph{polyads} in the chemistry literature. While we can
recognize the footprints of these polyads on the spectrum of
$\hat{P}_\varepsilon$, we notice that much of the structure is lost,
even for relatively small $\varepsilon$.

\subsection{$V_\varepsilon$ has a  global minimum}
Our goal is to compare the spectrum of $P_\varepsilon$ with what the
$\varepsilon$-Birkhoff-Gustavson normal form gives. In order to apply
our results, let us first check that $V_\varepsilon$ has a unique
non-degenerate minimum, when the constants $\gamma,c,\varepsilon$ are
properly chosen.

\begin{prop}\label{prop:unique-min}
  Assume that, in addition to~\eqref{equ:V-min-non-deg}, the following
  conditions hold:
  \begin{gather}
    1-\frac{\gamma^2}{8(1+c\varepsilon)} >0 \\
    \frac{7}{16}+\varepsilon-\frac{1}{4\left(1-\frac{\gamma^2}{8(1+c\varepsilon)}\right)}>0\,.
  \end{gather}
  Then, the potential $V_\varepsilon$ admits a unique non-degenerate
  minimum at the origin. In particular, if
  $\abs{\gamma}<\sqrt{\frac{24}{7}}$ and $\abs{\varepsilon}$ is small
  enough, then $V_\varepsilon$ admits a unique non-degenerate minimum
  at the origin.
\end{prop}
\begin{demo}
  We simply complete three squares: first, we write
  \[
    \tfrac{1}{2}x_1x_2^2 + x_4^4 = \left(x_2^2+\frac{x_1}{4}\right)^2
    - \frac{x_1^2}{16}\,.
  \]
  Then, we use
  \[
    (2+2c\varepsilon) x_2^2 + \gamma x_1^2 x_2 = 2(1+c\varepsilon)
    \left(x_2+\frac{\gamma x_1^2}{4(1+c\varepsilon)}\right)^2 -
    \frac{\gamma^2 x_1^4}{8(1+c\varepsilon)}.
  \]
  Incorporating the last term with the monomial $x_1^4$ from
  $V_\varepsilon$, we finally write, with
  $\kappa:=1-\frac{\gamma^2}{8(1+c\varepsilon)}$,
  \[
    \kappa x_1^4 + x_1^3 =
    x_1^2\left(\kappa\left(x_1+\frac{1}{2\kappa}\right)^2 -
      \frac{1}{4\kappa}\right)\,.
  \]
  This finally gives
  \begin{align}\label{equ:V-sum-of_squares}
    V_\varepsilon = & \left(\frac{1}{2} + \varepsilon - \frac{1}{16}
                      -\frac{1}{4\kappa}\right)x_1^2 + 2(1+c\varepsilon)
                      \left(x_2+\frac{\gamma x_1^2}{4(1+c\varepsilon)}\right)^2 \\
                    & + \kappa x_1^2\left(x_1+\frac{1}{2\kappa}\right)^2 +
                      \left(x_2^2+\frac{x_1}{4}\right)^2\,.
  \end{align}
  Thus, under the conditions of the proposition, we have a sum of four
  non-negative terms; the sum vanishes if and only if all terms
  vanish, which is equivalent to $x_1=x_2=0$.

  Now, if $\gamma$ and $c$ are fixed, the two conditions take the form
  \begin{gather}
    1-\frac{\gamma^2}{8} + \mathcal{O}(\varepsilon) >0\\
    \frac{7}{16} - \frac{1}{4(1-\frac{\gamma^2}{8})} +
    \mathcal{O}(\varepsilon) >0
  \end{gather}
  which holds, if $\varepsilon$ is small enough, as soon as
  \[
    \frac{7}{16} - \frac{1}{4(1-\frac{\gamma^2}{8})} > 0 \quad \text{
      i.e. } \quad \gamma^2 < \frac{24}{7}\,,
  \]
  which is stronger than the first condition $\gamma^2<8$.
\end{demo}
\begin{rema}
  We don't claim that the two conditions of
  Proposition~\ref{prop:unique-min} are necessary. But they allow for
  a simple proof, and they are sufficient for our numerical purposes
  $\gamma=0$ and $\gamma=1$.
\end{rema}

\subsection{Joint spectrum of $\hat{H}_{2,0}$ and $\hat{K}_{3,\varepsilon}$}
Let us apply the $\varepsilon$-Birkhoff-Gustavson normal form of
Theorem~\ref{theo:BGNF} to order 3: we conjugate the initial
Schrödinger operator $\hat P_\varepsilon$ to an operator of the form
$ H_{2,0}+K_{3,\varepsilon}+\mathcal{O}_4 $, and our goal is perform
numerical computations neglecting the $\mathcal{O}_4 $ term.

The operator $K_{3,\varepsilon}$ can be explicitly implemented thanks
to Theorem~\ref{theo:mu-reso12}, Theorem~\ref{theo:mu_explicit} and
Theorem~\ref{theo:matrix}. For any given $N\geq 0$, we obtain the
exact matrix for $K_{3,\varepsilon}$ in the orthonormal
basis~\eqref{equ:eigenspace}, of cardinal
$ \lfloor \tfrac{N}{2} \rfloor + 1$.

In order to obtain an approximation of the spectrum of
$\hat P_\varepsilon$, following Section~\ref{sec:spectr-analys-fermi},
we fix some energy $E>0$, and we compute the spectrum of
$H_{2,0}+K_{3,\varepsilon}$ restricted to the spectral subspace of
$H_{2,0}$ corresponding to energies in $[0,(1+\eta)E]$, where $\eta>0$
is fixed. In view of the Bargmann representation of
Section~\ref{sec:bargm-repr}, this is equivalent to computing the
spectrum of the restriction of
$\hat{H}_{2,0}^{\mathcal{B}} +\hat{K}_{3,\varepsilon}^{\mathcal{B}}$
to the space $\oplus_{n=0}^N\mathcal{H}_n^\mathcal{B}$,
see~\eqref{equ:eigenspace}, for some $N$ large enough.  Specifically,
the integer $N$ has to be chosen such that
$\h(\tilde
n+\frac{3}{2})+\sigma((\hat{K}_{3,\varepsilon}^{\mathcal{B}})_{\restr
  \mathcal{H}_{\tilde n}^\mathcal{B}})$ does not intersect the
interval $(0,E)$ for all $\tilde n > N$.

A nice way of displaying this computation is to plot the \emph{joint
  spectrum} of the commuting operators $\h^{-1}\hat{H}_{2,0}$ and
$\hat{K}_{3,\varepsilon}$: this is the set of pairs
$(\lambda_1,\lambda_2)\in\RM^2$ such that, with the notation of
Theorem~\ref{theo:matrix},
\begin{align}
  \lambda_1 & = \omega_{1,0}(N+\tfrac{3}{2})\\
  \lambda_2 & = \mu_{\varepsilon,N,j} \text{ for some } j=0,\dots,  \lfloor \tfrac{N}{2}
              \rfloor\,,
\end{align}
see Figure~\ref{fig:jspec}.  Then, the spectrum
$H_{2,0}+K_{3,\varepsilon}$ is obtained from the joint spectrum by
applying the map $(\lambda_1,\lambda_2)\mapsto \h\lambda_1+\lambda_2$.

Notice that the computation of the joint spectrum is much faster than
the computation of the spectrum of $\hat{P}_\varepsilon$ from
Section~\ref{sec:numer-comp-spectr}, since instead of a matrix of size
at least $\mathcal{O}(\frac{E^2}{\h^2}\times\frac{E^2}{\h^2} )$, we
have $N=\mathcal{O}({\frac{E}{\h}})$ matrices of sizes
$\mathcal{O}(\frac{E}{\h}\times\frac{E}{\h} )$, so we gain (at least)
an order of magnitude in $\frac{E}{\h}$. For instance, with $M=41$,
instead of a $903\times 903$ matrix, we have 10 matrices of sizes
$1\times 1, 2\times 2, 2\times 2, \dots, 6\times 6$.

\subsection{The case $\gamma=0$}
\label{sec:case-gamma=0}

The case $\gamma=0$ is interesting because the remaining terms of
order 3 in the potential $V_\varepsilon$, namely $x_1^3$ and
$x_1x_2^2$, are completely cancelled out by the Birkhoff normal form:
we have $\mu=0$ in Theorem~\ref{theo:mu_explicit}. Therefore, the
spectrum of $\hat{P}_\varepsilon$ is approximated up to
$\mathcal{O}_4$ merely by the quadratic term
$\hat{H}_{2,0}+\varepsilon\hat{L}_2$, whose spectrum is explicit:
\[
  \left\{\h\left(\omega_{1,0}(\alpha_1 + 2\alpha_2 +
      \tfrac{3}{2})\right) +
    \h\varepsilon\left(\omega_{1,1}(\alpha_1+\tfrac{1}{2}) +
      \omega_{2,1}(\alpha_2+\tfrac{1}{2})\right),\quad
    (\alpha_1,\alpha_2)\in \NM^2\right\}
\]

\begin{figure}[h]
  \centering \includegraphics[width=0.45\textwidth]{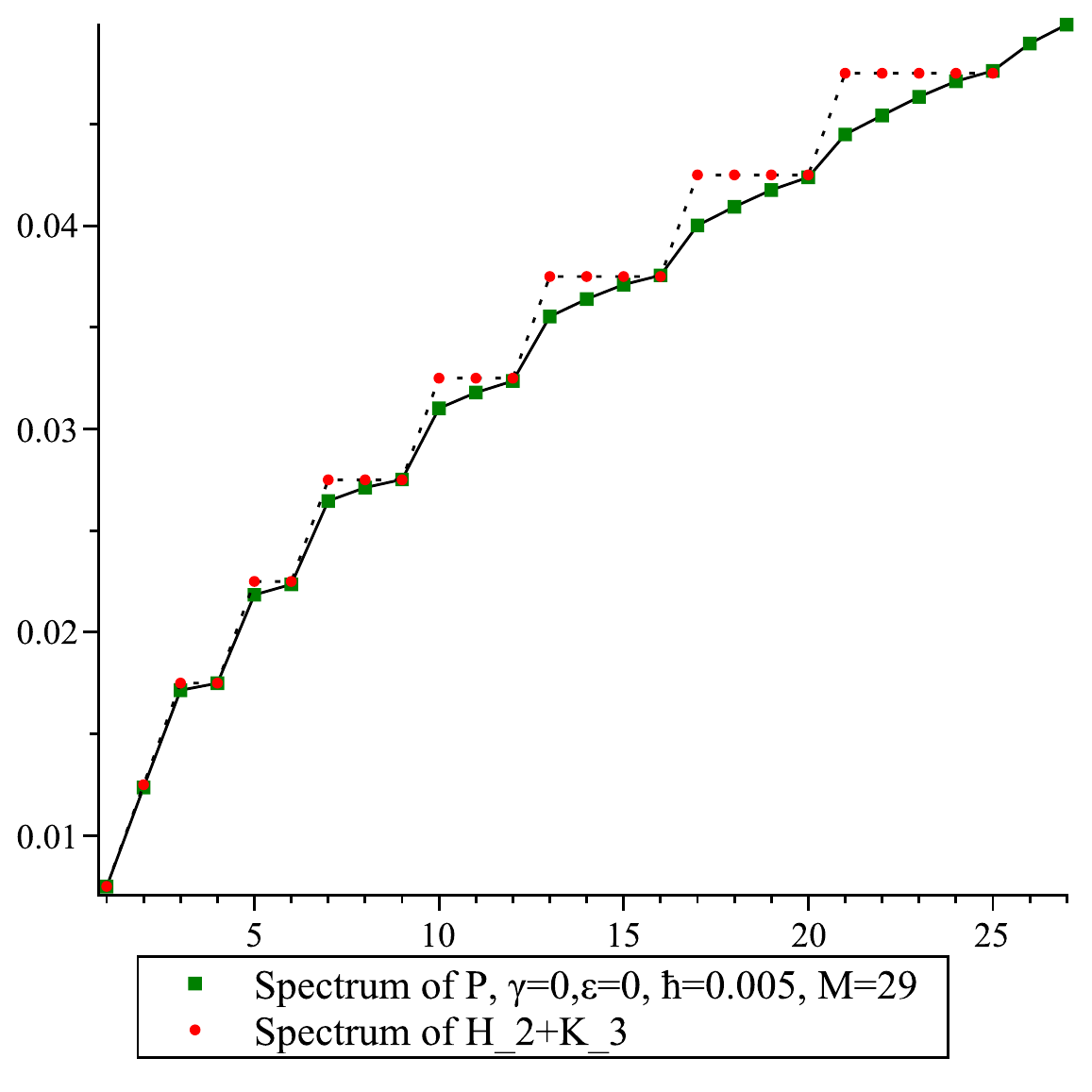}
  \includegraphics[width=0.45\textwidth]{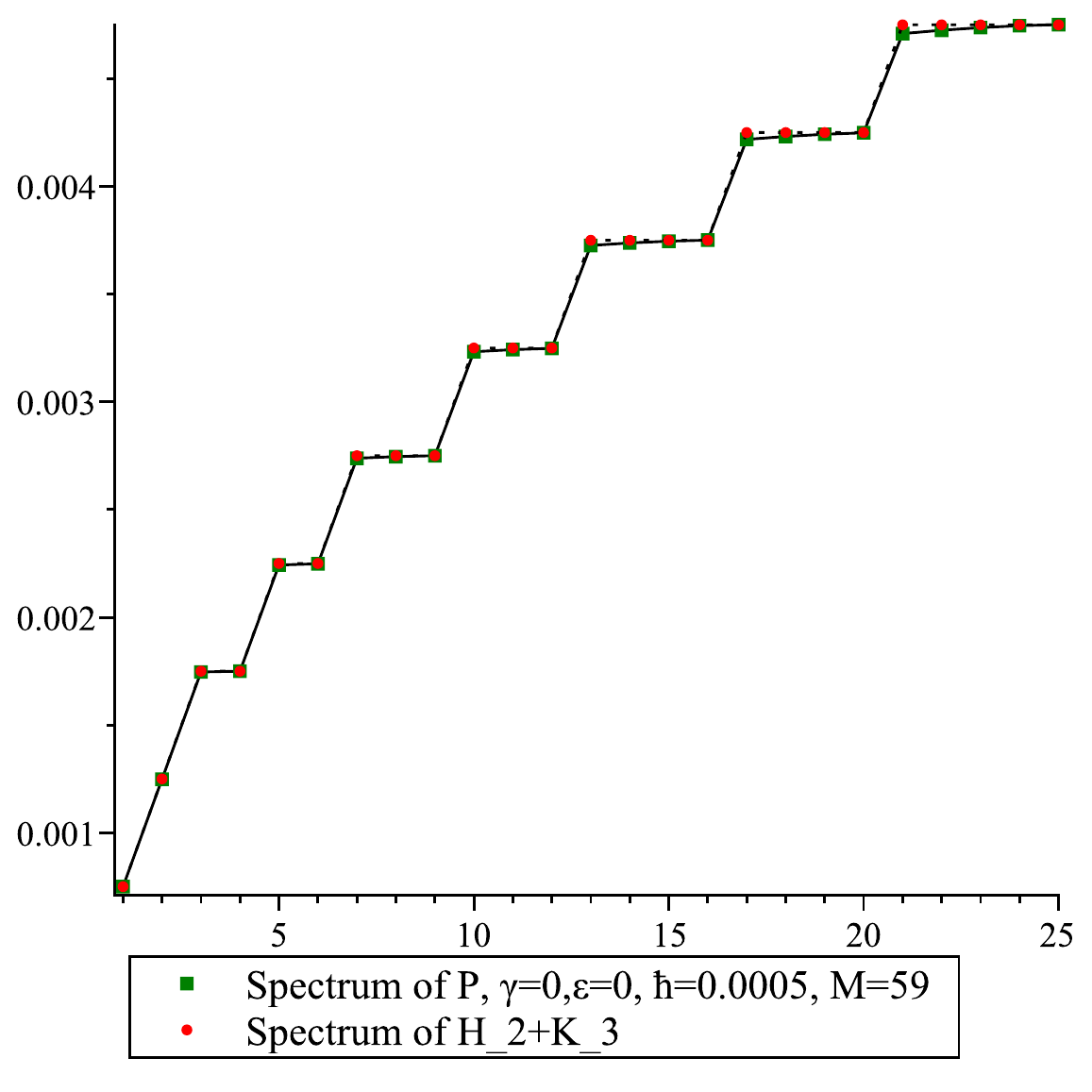}
  \caption{Comparison of the spectrum of $\hat{P}_\varepsilon$ with
    the eigenvalues obtained from the Birkhoff normal form
    $ H_{2,0}+K_{3,\varepsilon}$. Here $\h=0.005$ (left) or
    $\h=0.0005$ (right), and $\gamma=0,\varepsilon=0$, hence
    $K_{3,\varepsilon}=0$.}
  \label{fig:g=0_e=0}
\end{figure}

\begin{figure}[h]
  \centering \includegraphics[width=0.5\textwidth]{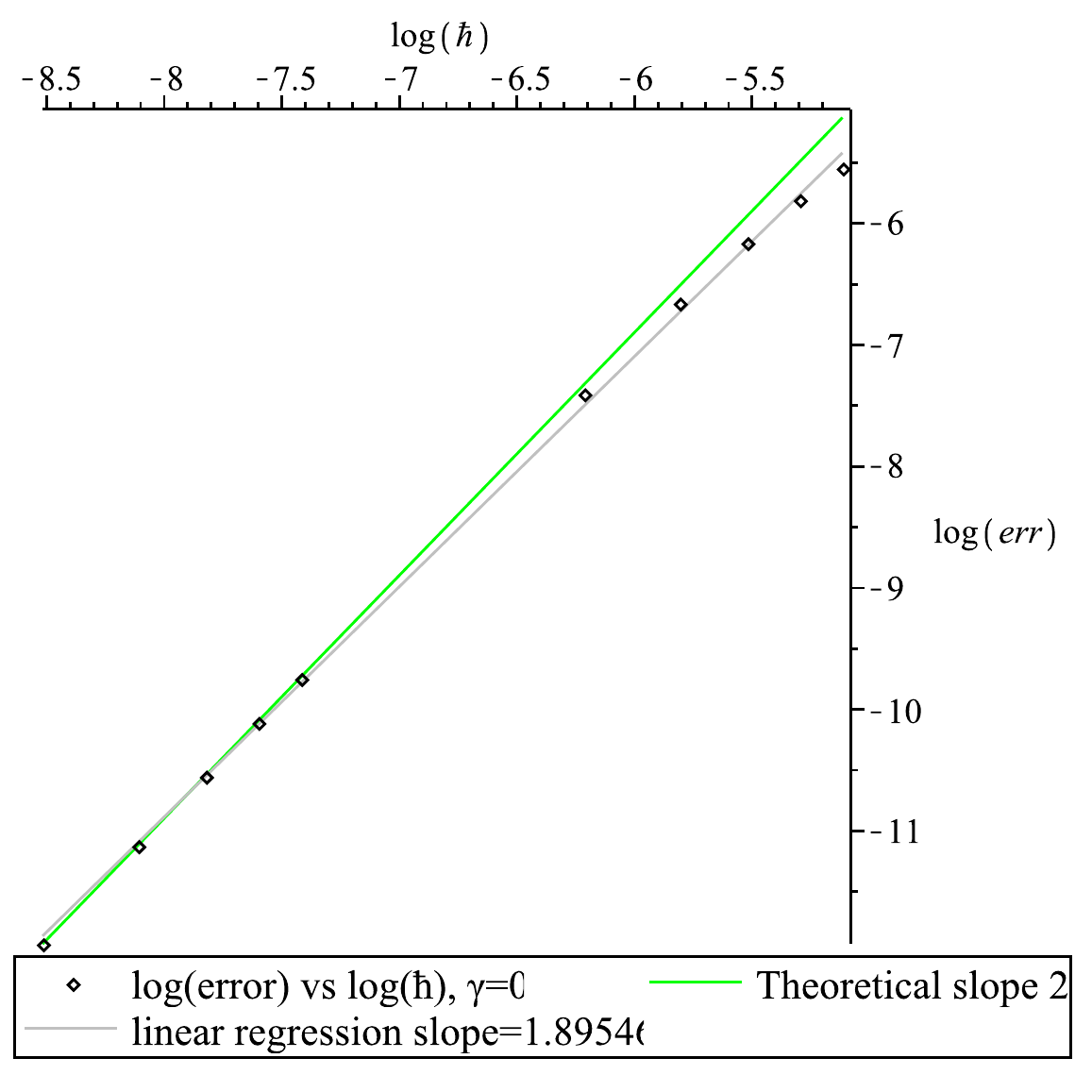}
  \caption{Error(log scale) between the spectrum of
    $\hat{P}_\varepsilon$ and that of $
    H_{2,0}+K_{3,\varepsilon}$. Here $\gamma=0,\varepsilon=0$. When
    $\h$ is small enough, we observe the theoretical slope of 2 (green
    line), corresponding to an error of order $\mathcal{O}(\h^2)$.}
  \label{fig:g=0_e=0_err}
\end{figure}

When $\varepsilon=0$, the theoretical results of~\cite{san-charles}
apply. In particular, in the regime $E=C\h$, the spectrum should
converge to the polyads of~\eqref{equ:polyads}, with an error of order
$\mathcal{O}(\h^2)$. This is clearly illustrated in
Figures~\ref{fig:g=0_e=0} and~\ref{fig:g=0_e=0_err}.

\begin{figure}[h]
  \centering \includegraphics[width=0.45\textwidth]{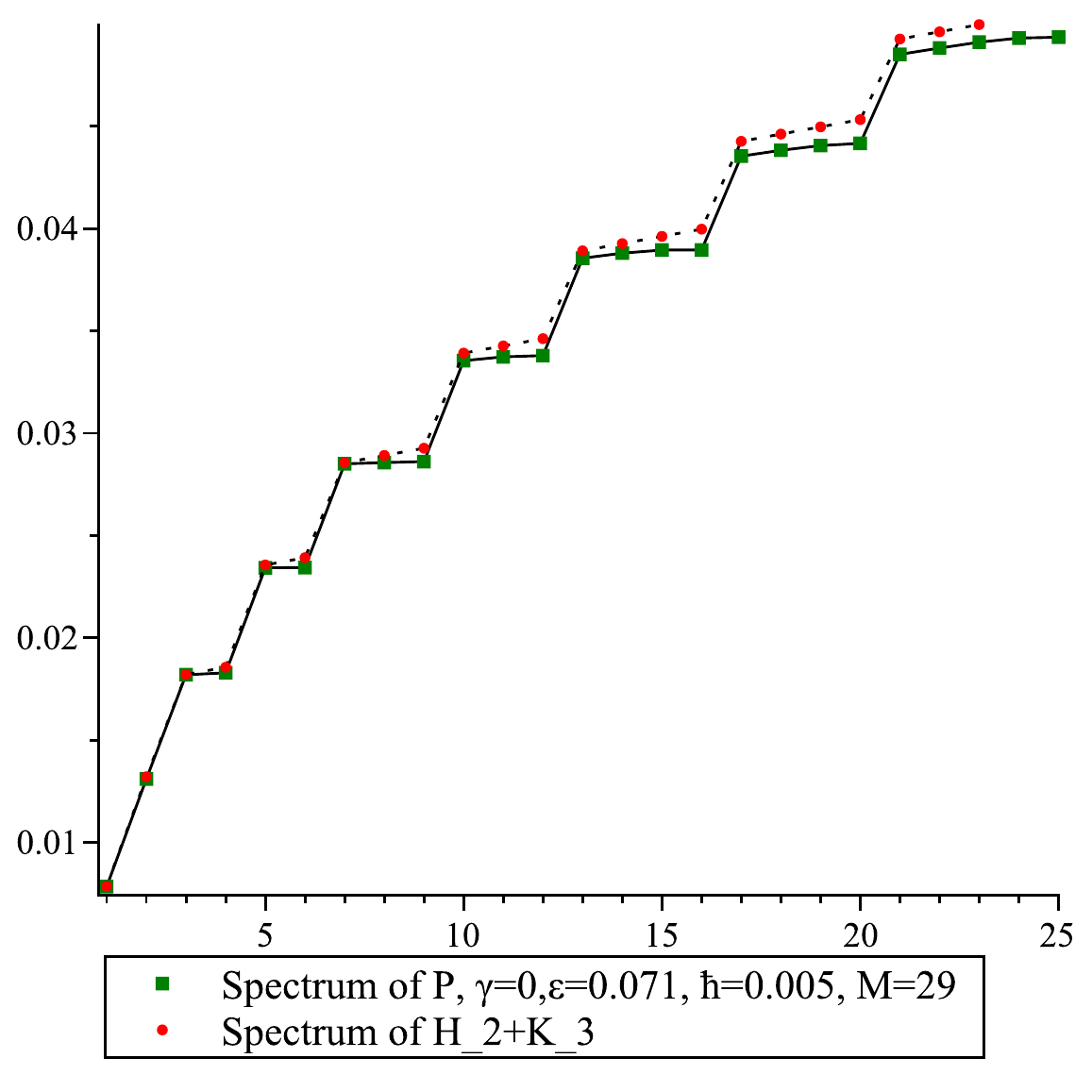}
  \includegraphics[width=0.45\textwidth]{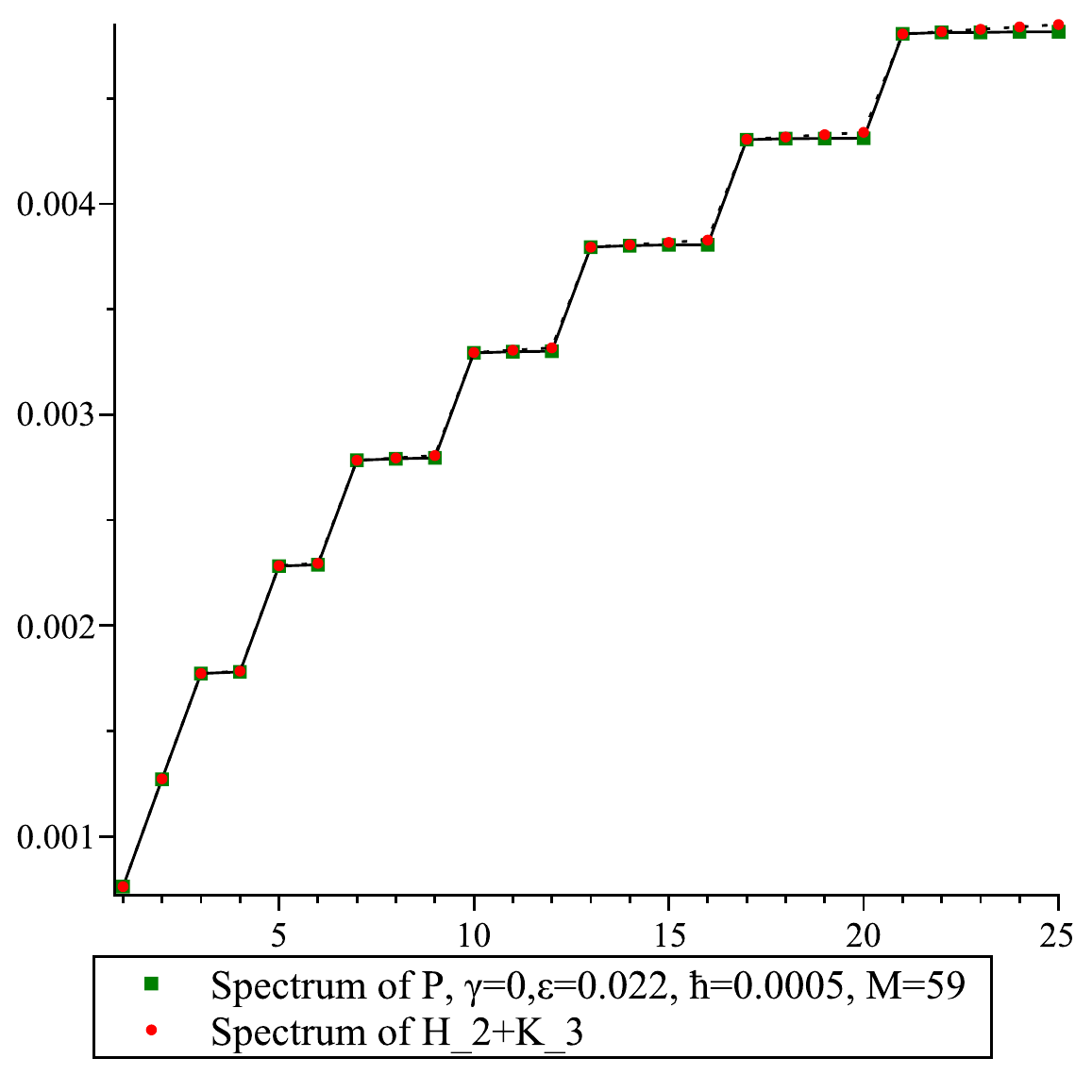}
  \caption{Comparison of the spectrum of $\hat{P}_\varepsilon$ with
    the eigenvalues obtained from the Birkhoff normal form
    $ H_{2,0}+K_{3,\varepsilon}$. Here $\h=0.005$ (left) or
    $\h=0.0005$ (right), and $\gamma=0,\varepsilon=\sqrt{\h}$, hence
    $K_{3,\varepsilon}=0$.}
  \label{fig:g=0_e=sqrt}
\end{figure}

\begin{figure}[h]
  \centering \includegraphics[width=0.5\textwidth]{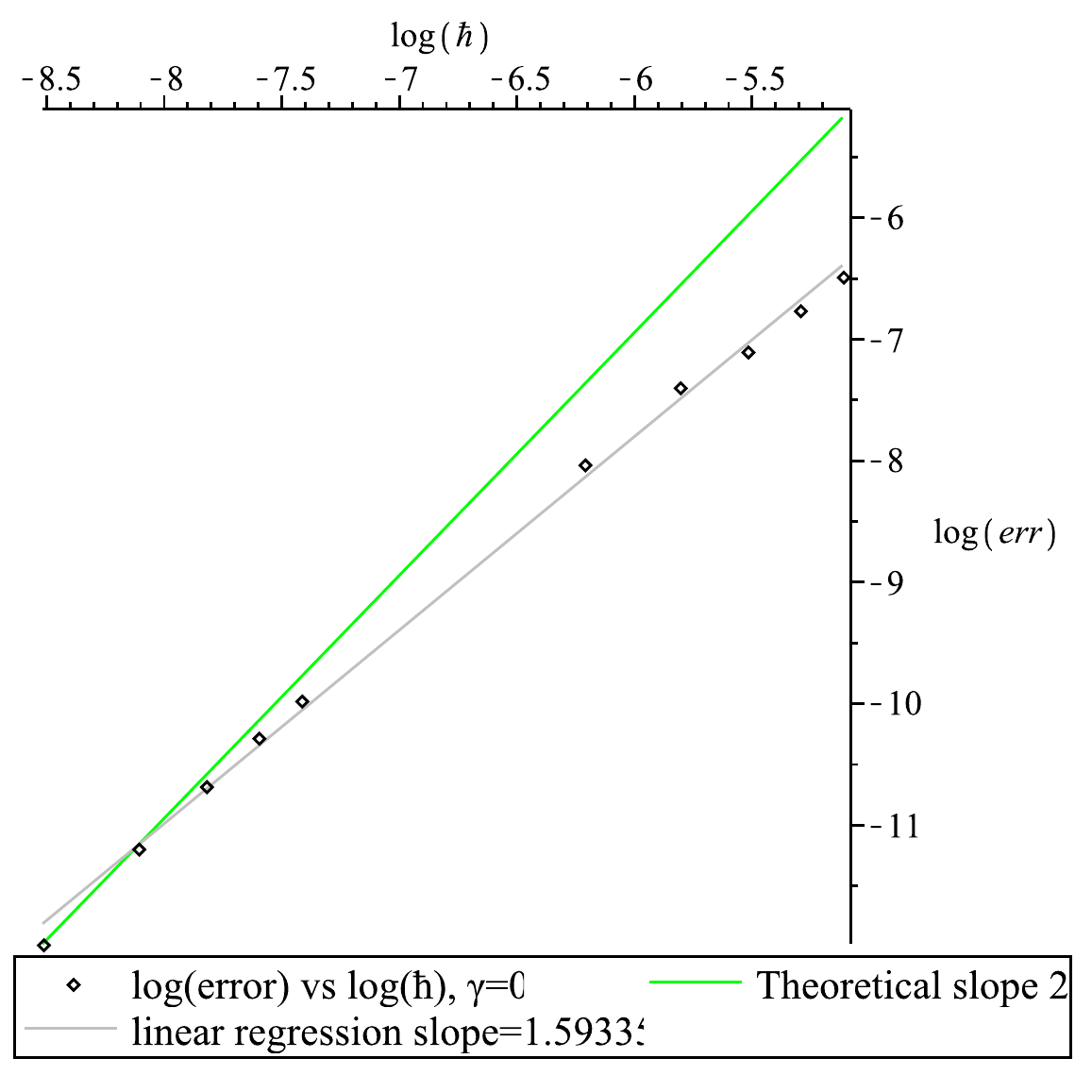}
  \caption{Error(log scale) between the spectrum of
    $\hat{P}_\varepsilon$ and that of $
    H_{2,0}+K_{3,\varepsilon}$. Here $\gamma=0,\varepsilon=\sqrt{h}$.
    When $\h$ is small enough, we tend to the theoretical slope of 2
    (green line), corresponding to an error of order
    $\mathcal{O}(\h^2)$, though it seems slower than in the case
    $\varepsilon=0$ (Figure~\ref{fig:g=0_e=0_err}); on the other hand,
    for the larger $\h$, the absolute error is quite smaller here than
    for the case $\varepsilon=0$.}
  \label{fig:g=0_e=sqrt_err}
\end{figure}

In order to experiment the case $\varepsilon\neq 0$, we chose the
regime $\varepsilon=\sqrt{\h}$, where we still expect an error of
order $\mathcal{O}(\h^2)$, which is confirmed by the numerics, see
Figures~\ref{fig:g=0_e=sqrt} and~\ref{fig:g=0_e=sqrt_err}.

\begin{figure}[h]
  \centering
  \includegraphics[width=0.7\linewidth]{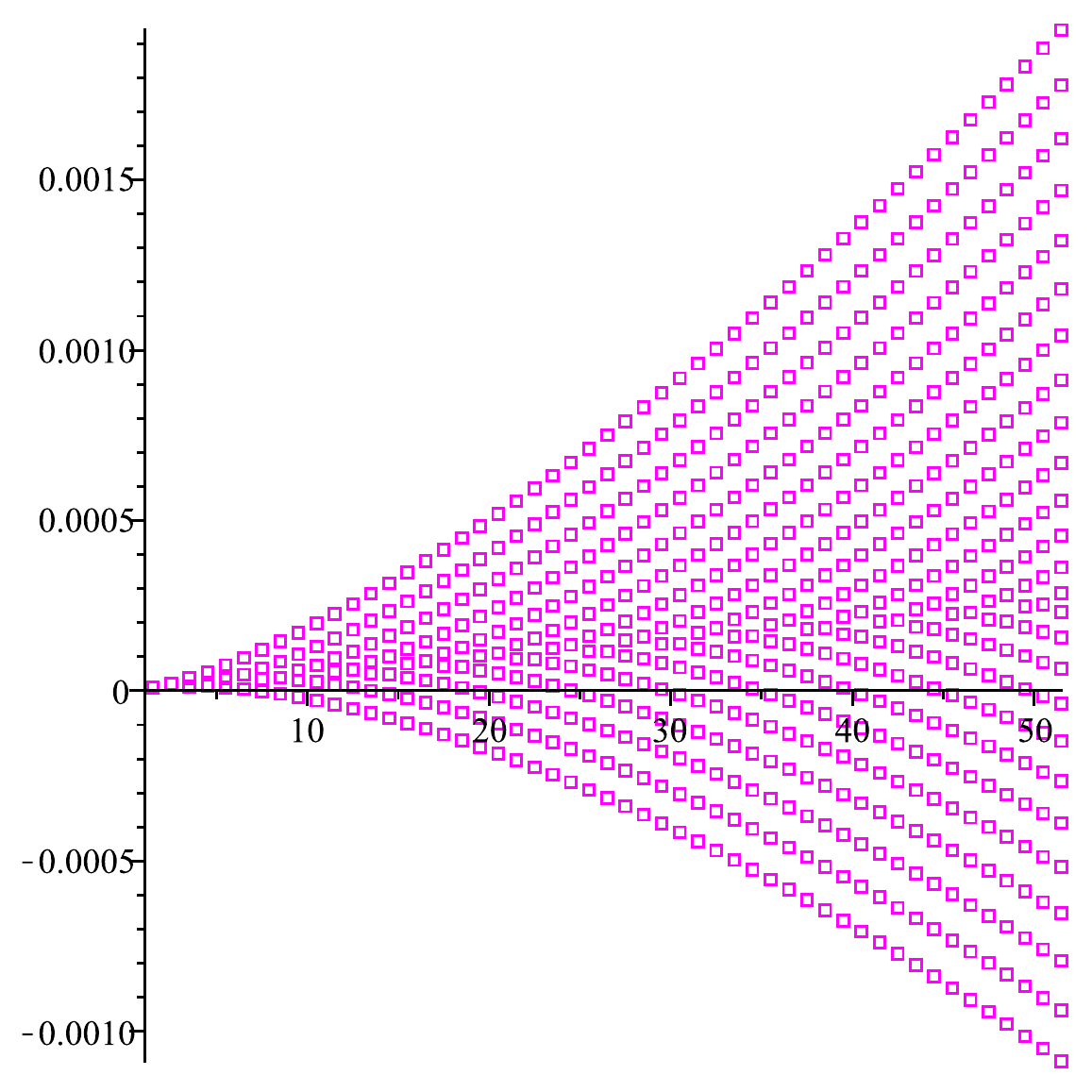}
  \caption{Joint spectrum of
    $(\h^{-1}\hat{H}_{2,0} , \hat{K}_{3,\varepsilon})$. The abscissae
    are $n+\frac{3}{2}$, $n\geq 0$. Here $\gamma=1,c=1$,
    $\varepsilon=0.01, \h=0.001$. }
  \label{fig:jspec}
\end{figure}

The joint spectrum of the commuting operators
$(\h^{-1}\hat{H}_{2,0} , \hat{K}_{3,\varepsilon})$ is depicted in
Figure~\ref{fig:jspec}.

\subsection{The case $\gamma=1$}

The case $\gamma=1$ corresponds to the heart of our result, since the
Birkhoff term of order 3, $K_{3,\varepsilon}$ is not trivial, due
to~\eqref{equ:mu-exemple}.

\begin{figure}[h]
  \centering \includegraphics[width=0.45\textwidth]{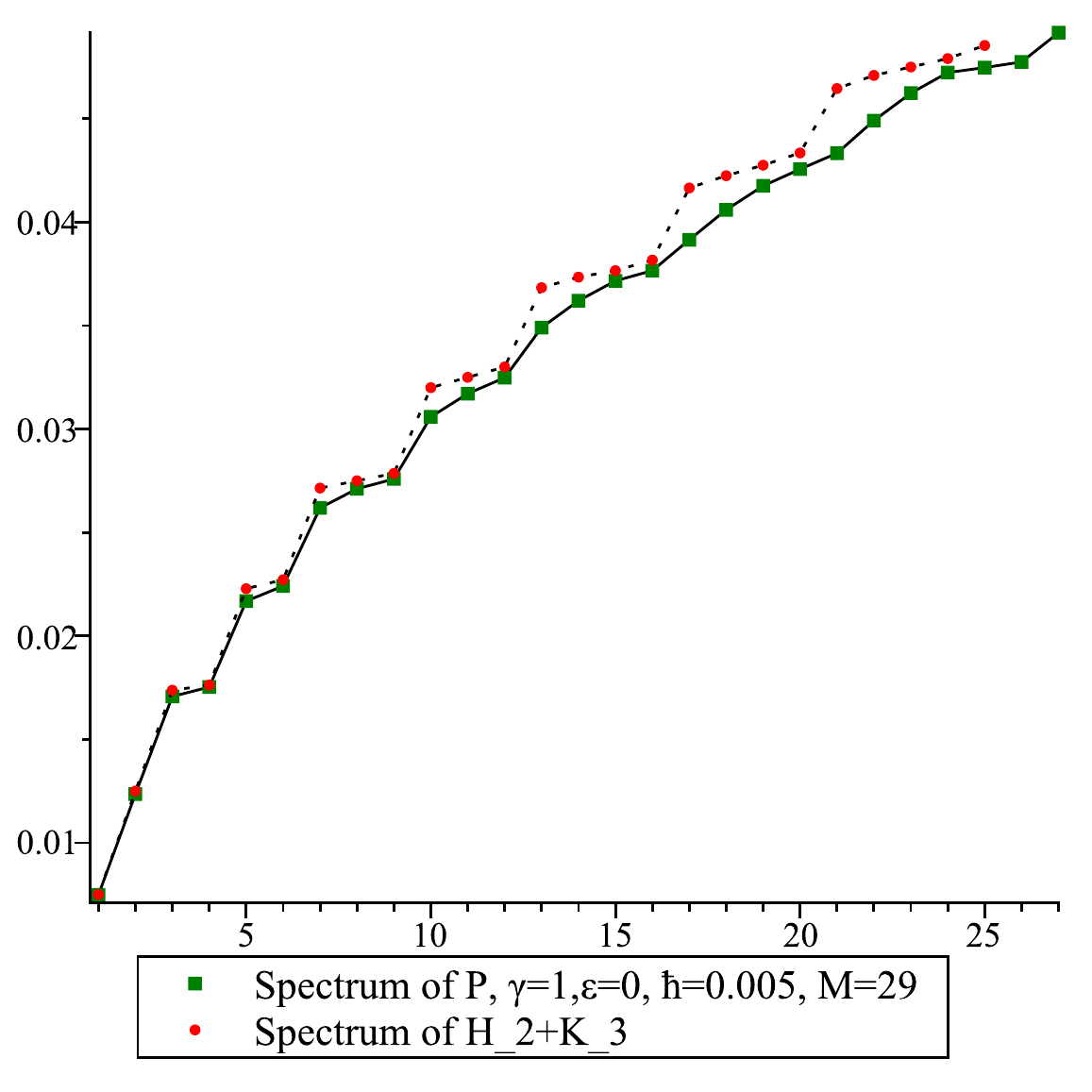}
  \includegraphics[width=0.45\textwidth]{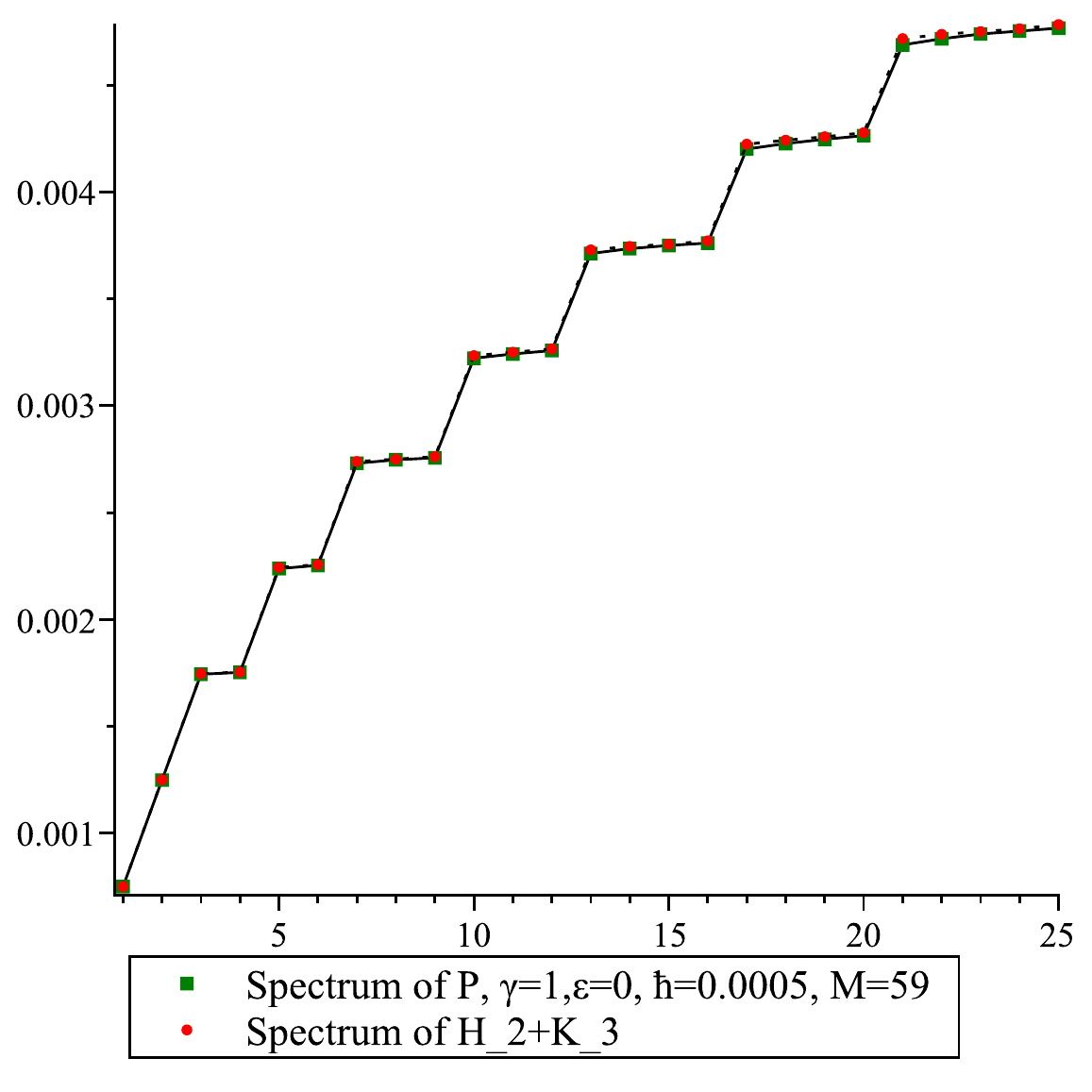}
  \caption{Comparison of the spectrum of $\hat{P}_\varepsilon$ with
    the eigenvalues obtained from the Birkhoff normal form
    $ H_{2,0}+K_{3,\varepsilon}$. Here $\h=0.005$ (left) or
    $\h=0.0005$ (right), and $\gamma=1,\varepsilon=0$.}
  \label{fig:g=1_e=0}
\end{figure}

\begin{figure}[h]
  \centering \includegraphics[width=0.5\textwidth]{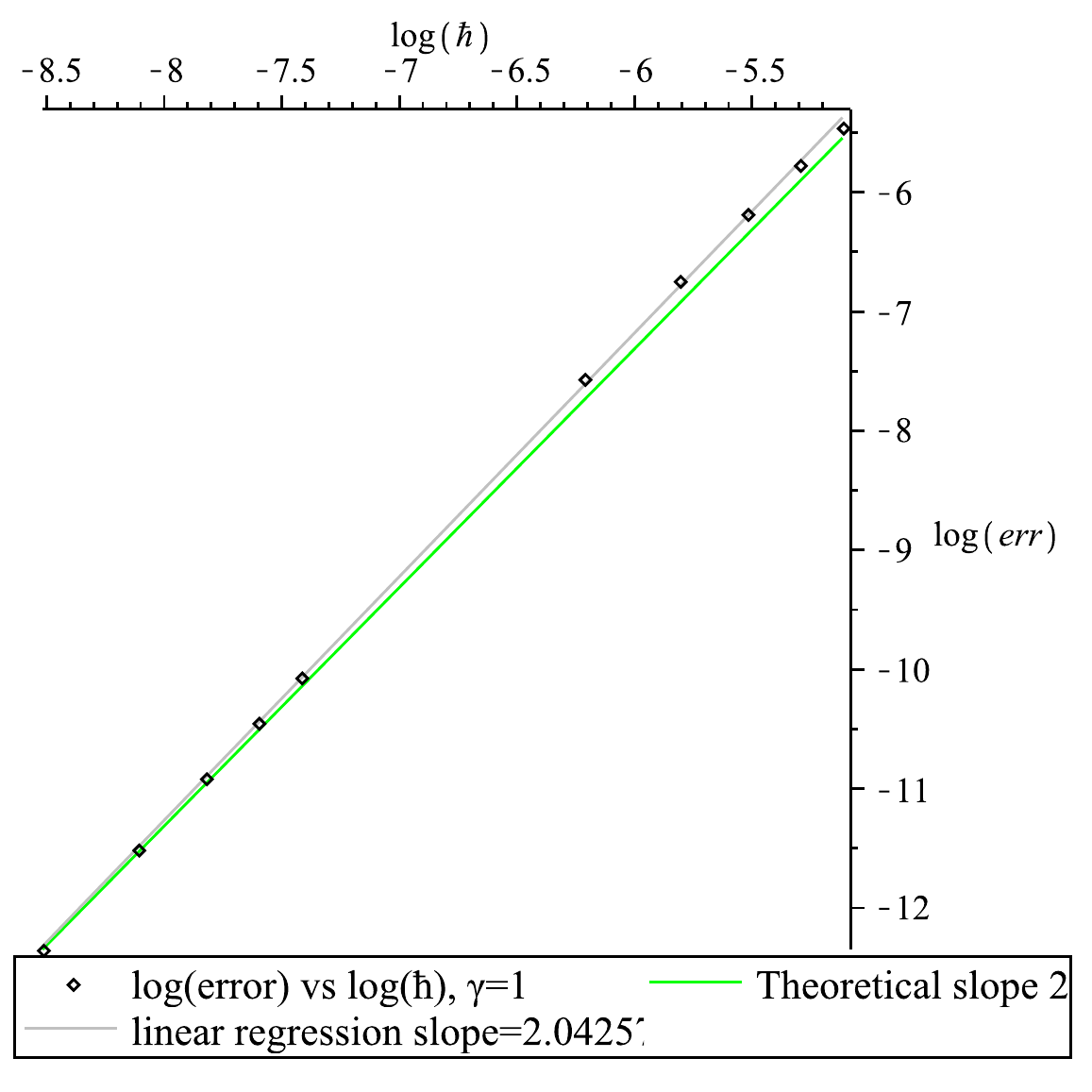}
  \caption{Error(log scale) between the spectrum of
    $\hat{P}_\varepsilon$ and that of $
    H_{2,0}+K_{3,\varepsilon}$. Here $\gamma=1,\varepsilon=0$. We
    observe a nice fit with the theoretical slope of 2 (green line),
    corresponding to an error of order $\mathcal{O}(\h^2)$.}
  \label{fig:g=1_e=0_err}
\end{figure}

When $\varepsilon=0$, as above, the theoretical results
of~\cite{san-charles} apply and we observe the expected
$\mathcal{O}(\h^2)$ error in Figures~\ref{fig:g=1_e=0}
and~\ref{fig:g=1_e=0_err}. The clear agreement with
$\mathcal{O}(\h^2)$ is a strong confirmation of the validity of the
Birkhoff procedure, and in particular of the correctness of the value
of $\mu$ from~\eqref{equ:mu-exemple}, because any other value of $\mu$
would lead to an error of the order of the eigenvalues of $K_{3,0}$ on
the given spectral subspace, which is known to be
$\mathcal{O}(\h^{3/2})$.

\begin{figure}[h]
  \centering \includegraphics[width=0.45\textwidth]{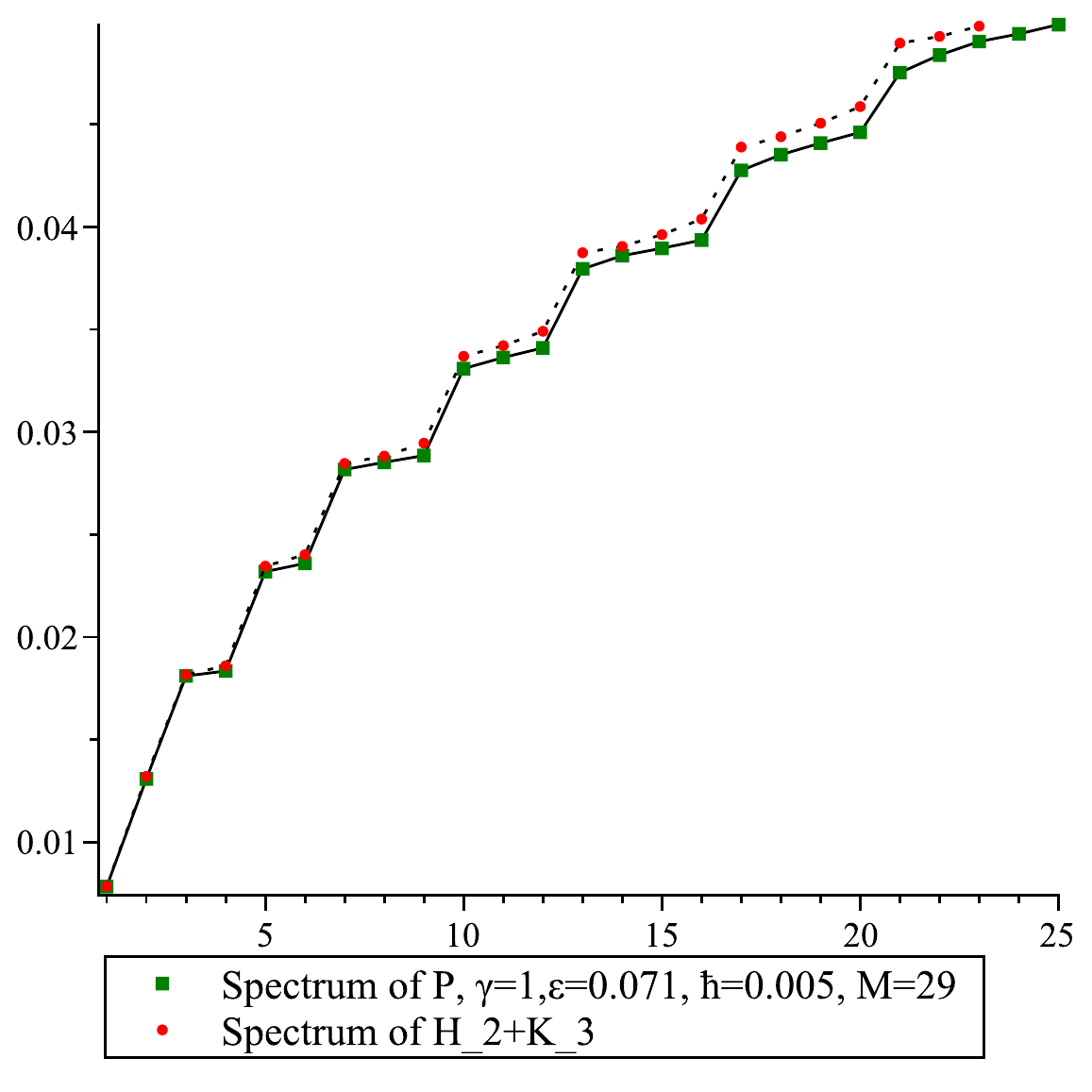}
  \includegraphics[width=0.45\textwidth]{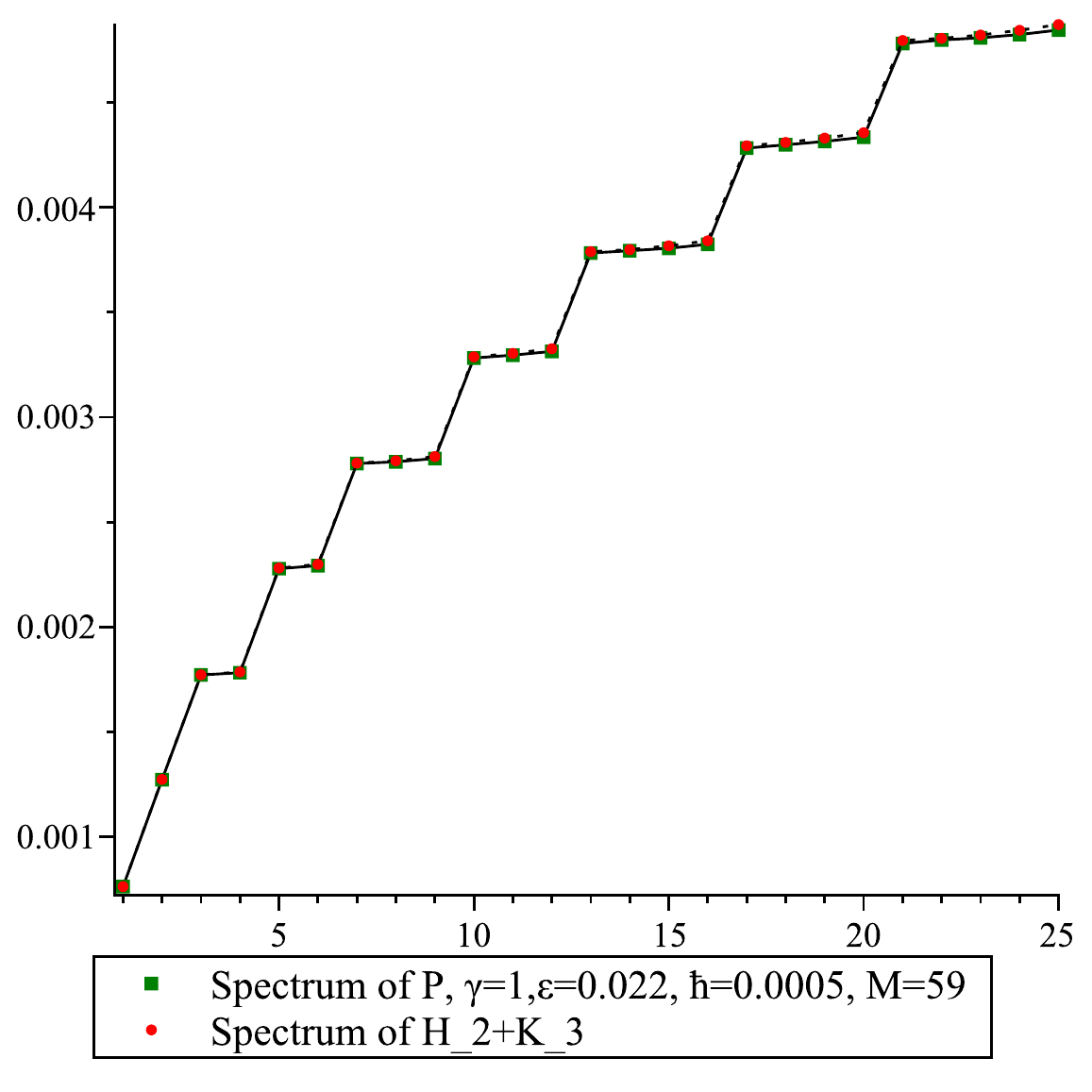}
  \caption{Comparison of the spectrum of $\hat{P}_\varepsilon$ with
    the eigenvalues obtained from the Birkhoff normal form
    $ H_{2,0}+K_{3,\varepsilon}$. Here $\h=0.005$ (left) or
    $\h=0.0005$ (right), and $\gamma=1,\varepsilon=\sqrt{\h}$.}
  \label{fig:g=1_e=sqrt}
\end{figure}

\begin{figure}[h]
  \centering \includegraphics[width=0.5\textwidth]{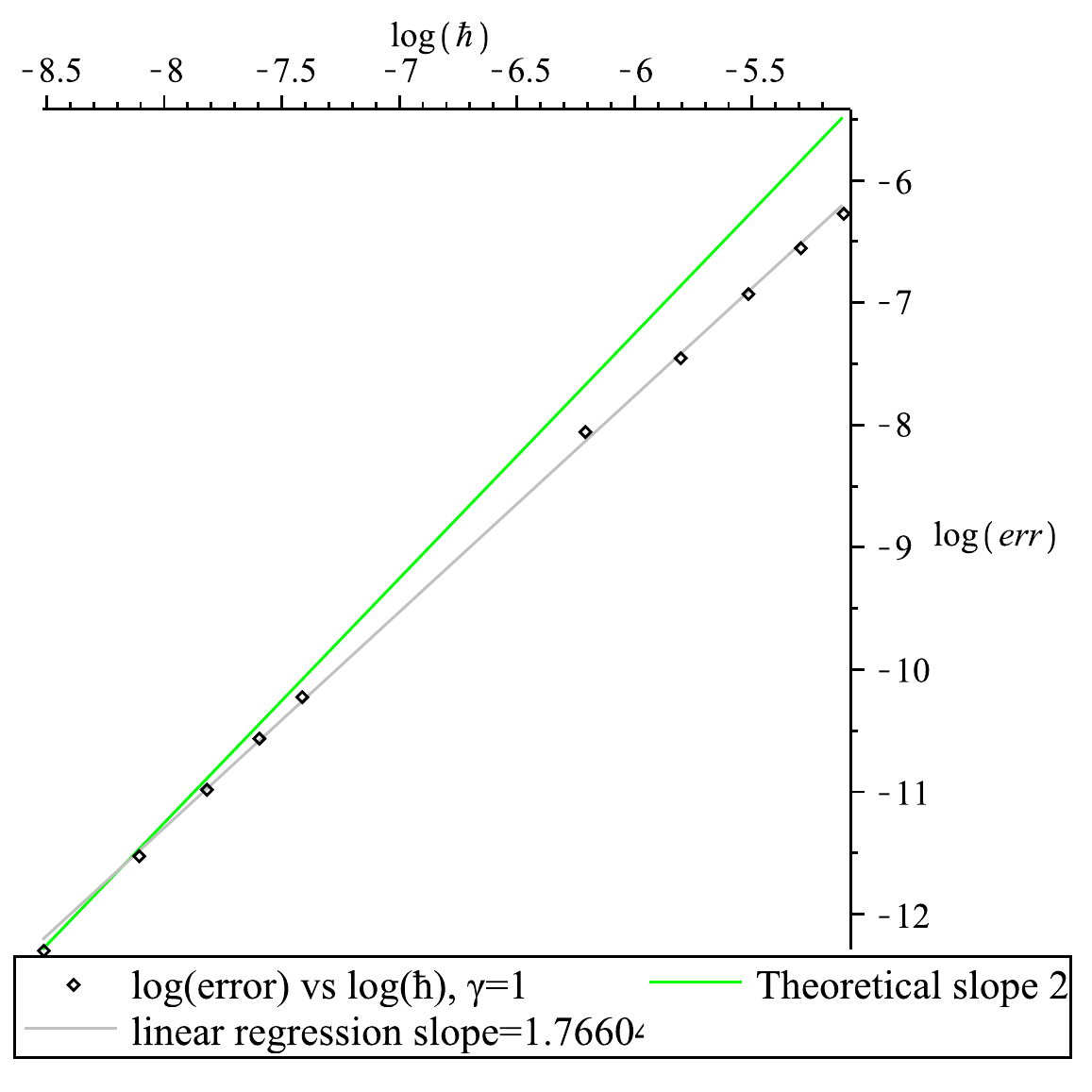}
  \caption{Error(log scale) between the spectrum of
    $\hat{P}_\varepsilon$ and that of $
    H_{2,0}+K_{3,\varepsilon}$. Here
    $\gamma=1,\varepsilon=\sqrt{\h}$. When $\h$ is small enough, we
    tend to the theoretical slope of 2 (green line), corresponding to
    an error of order $\mathcal{O}(\h^2)$, though it seems slower than
    in the case $\varepsilon=0$ (Figure~\ref{fig:g=1_e=0_err}); on the
    other hand, for the larger $\h$, the absolute error is quite
    smaller here than for the case $\varepsilon=0$.}
  \label{fig:g=1_e=sqrt_err}
\end{figure}

The new results correspond to $\varepsilon\neq 0$. As in the case
$\gamma=0$ we experiment the regime $\varepsilon=\sqrt{\h}$, and, in
spite of this perturbation, the $\varepsilon$-Birkhoff-Gustavson
procedure suggests that the error should still be
$\mathcal{O}(\h^2)$. This is confirmed by Figures~\ref{fig:g=1_e=sqrt}
and~\ref{fig:g=1_e=sqrt_err}.

\clearpage

\paragraph{Acknowledgements.}
This work is supported in part by funds provided by Henri Lebesgue
Center (ANR LEBESGUE), the Laboratory of Fundamental and Applied
Mathematics of Oran and the Algerian research project: PRFU No:
C00L03ES310120220001. The second author is happy to acknowledge the
excellent working conditions that she was given during several stays
at the IRMAR.

\bibliographystyle{abbrv}%
\bibliography{bibli-utf8}
\end{document}